\documentclass[10pt]{article}
\usepackage{latexsym,amsmath,amssymb,amscd, graphics}

\newcommand{\bfi}{\bfseries\itshape}

\usepackage{amsmath}

\makeatletter
\@addtoreset{figure}{section}
\def\thefigure{\thesection.\@arabic\c@figure}
\def\fps@figure{h,t}
\@addtoreset{table}{bsection}

\def\thetable{\thesection.\@arabic\c@table}
\def\fps@table{h, t}
\@addtoreset{equation}{section}

\makeatother

\newtheorem{theorem}{Theorem}

\newtheorem{corollary}[theorem]{Corollary}

\newtheorem{definition}[theorem]{Definition}
\newtheorem{example}[theorem]{Example}

\newtheorem{lemma}[theorem]{Lemma}

\newtheorem{proposition}[theorem]{Proposition}
\newtheorem{remark}[theorem]{Remark}

\numberwithin{theorem}{section}
\numberwithin{equation}{section}

\begin{document}

\title{Openness and convexity for momentum maps}
\author{Petre Birtea, Juan-Pablo Ortega, Tudor S. Ratiu}
\date{}
\maketitle

\begin{abstract}
The purpose of this paper is finding the essential attributes
underlying the convexity theorems for momentum maps. It is shown
that they are of topological nature; more specifically, we show
that convexity follows if the map is open onto its image and has
the so called local convexity data property. These conditions are
satisfied in all the classical convexity theorems and hence they
can, in principle, be obtained as corollaries of a more general
theorem that has only these two hypotheses. We also prove a
generalization of the ``Lokal-global-Prinzip" that only requires
the map to be closed and to have a normal topological space as
domain, instead of using a properness condition. This allows us to
generalize the  Flaschka-Ratiu convexity theorem to non-compact
manifolds.
\end{abstract}

\section{Introduction}

The problem of describing the image of a momentum map
defined on a symplectic manifold has generated a large amount of
research in the past twenty years and it remains to this day one
of the most active areas in symplectic geometry and
its applications to Hamiltonian dynamics, especially bifurcation
theory. In 1982 Atiyah \cite{atiyah 82} and, independently,
Guillemin and Sternberg \cite{convexity} proved the
following result about the convexity of the image of the
momentum map associated to the action of a torus $T$ on a
compact symplectic manifold.

\medskip

\noindent {\bfi Theorem (Atiyah-Guillemin-Sternberg)} \textit{
Let $M$ be a compact connected symplectic manifold on which a
torus $T$ acts in a Hamiltonian fashion with associated
invariant momentum map $J:M\rightarrow \mathfrak{t}^{*}$. Here
$\mathfrak{t}$ denotes the Lie algebra of $T$ and
$\mathfrak{t}^{\ast}$ is its dual; both are isomorphic as Abelian
Lie algebras to $\mathbb{R}^{{\rm dim} T}$. Then the image
$J(M)$ of $J$ is a compact convex polytope in ${\Bbb R}^{{\rm
dim} T}$, called the $T$-momentum polytope. Moreover, it is equal
to the convex hull of the image of the fixed point set of the
$T$-action. The fibers of $J$ are connected.}

\medskip

The motivation for this result was the \textit{Kostant linear
convexity theorem} \cite{kostant73} which states that the
projection of an adjoint orbit of a compact connected Lie group
onto the Lie algebra of any of its maximal tori is the convex hull
of the corresponding Weyl group orbit. This in turn is a
generalization of a classical result of Schur \cite{schur} and
Horn \cite{horn} (which in the case of the unitary group $U(n)$ is
Kostant's theorem) which states that the set of diagonals of an
isospectral set of Hermitian matrices equals the convex hull of
the $n!$ points obtained by permuting all the eigenvalues.

Another problem with important mathematical and physical
consequences is to describe all possible eigenvalues of the sum
$A + B $ of two Hermitian matrices $A $ and $B $ as each one of
them ranges over an isospectral set. The isospectral sets of
Hermitian matrices are precisely the coadjoint orbits
$\mathcal{O}_\mu$ of $U(n)$, where $\mu = (\mu_1, \mu_2, \dots,
\mu_n)$ and $\mu_i $ are the eigenvalues.  If one requires, in
addition, that the eigenvalues of the sum $A + B $ be sorted in
decreasing order, this problem amounts to describing the set
$J(\mathcal{O}_{\mu}\times \mathcal{O}_{\lambda })\cap
\mathfrak{t}_{+}^{*}$ where $J(A, B) = A +B$ and
$\mathfrak{t}_{+}^{*}$ is a positive Weyl chamber in
$\mathfrak{t}^\ast = \{\text{real diagonal matrices}\}$. This
problem is a particular case of the following more general
situation: a compact Lie group $G$ acts on a compact symplectic
manifold $M$ with associated equivariant momentum map $J:M\rightarrow
\mathfrak{g}^{*}$. Equivariance is not an assumption, since any
momentum map of a compact group can be averaged to give an
equivariant momentum map for the same canonical $G $-action.

Guillemin and Sternberg \cite{convexity2} proved that
$J(M)\cap \mathfrak{t}^{*}$ is a union of compact convex
polytopes and Kirwan \cite{kirwan convexity} showed that this set
is connected thereby concluding that $J(M)\cap
\mathfrak{t}_{+}^{*}$ is a compact convex polytope. We will
refer to this as the $G$-momentum polytope.
\medskip

\noindent {\bfi Theorem (Guillemin-Kirwan-Sternberg)} \textit{Let
$M$ be a compact connected symplectic manifold on which the
compact connected Lie group $G$ acts in a Hamiltonian fashion with
associated equivariant momentum map $J:M\rightarrow
\mathfrak{g}^{*}$. Here $\mathfrak{g}$ denotes the Lie algebra of
$G $ and $\mathfrak{g}^\ast$ is its dual. Let $T $ be a maximal
torus of $G $, $\mathfrak{t}$ its Lie algebra, $\mathfrak{t}^\ast$
its dual, and $\mathfrak{t}^\ast_+ $ the positive Weyl chamber
relative to a fixed ordering of the roots. Then $J(M)\cap
\mathfrak{t}^{\ast}_{+}$ is a compact convex polytope, called the
$G$-momentum polytope. The fibers of $J$ are connected.}
\medskip

Important results about the description of the $U(n)$-momentum
polytope were obtained by Knutson and Tao \cite{knotsontao}.

Sjamaar \cite{sjamaar} has given another proof of the
convexity theorems, based on ideas coming from K\"ahler and
algebraic geometry. His proof, even though it gives the most
complete information on these polytopes,  uses strongly the
symplectic form and it is not known how to generalize his
technique to other manifolds and other types of actions, such as
Poisson actions of Poisson-Lie groups.

The case of compact symplectic manifolds is rich but quite
particular. For non-compact manifolds, the previous results no
longer hold and a counterexample was given by Prato
\cite{prato}. Conditions under which the $T$ or $G$-momentum
polytopes are convex were given by Condevaux, Dazord, and Molino
\cite{dazord} and later by Hilgert, Neeb, and Plank
\cite{hilnebplank}. These papers show that the proof of the
convexity of the image of the momentum map rests on the
following result which we give here in the formulation due to the
second group of authors:

\medskip

\noindent {\bfi Theorem (Lokal-global-Prinzip)} \textit{ Let $\Psi
:X\rightarrow V$ be a locally fiber connected map from a connected
locally connected Hausdorff topological space $X$ to a finite
dimensional vector space $V$, with local convexity data
$(C_{x})_{x\in X}$ such that all convex cones $C_x$ are closed in
$V$. Suppose that $\Psi$ is a proper map. Then $\Psi(X)$ is a
closed locally polyhedral convex subset of $V$, the fibers $\Psi
^{-1}(v)$ are all connected, and $\Psi :X\rightarrow \Psi (X)$ is
an open mapping.}

\medskip

We elaborate now on the hypotheses of this theorem. A map $\Psi
:X\rightarrow V$ is said to have {\bfi local convexity data}
if for each $x\in X$ there exists an arbitrarily small open neighborhood $U_{x}$ of $x$ and a convex cone $C_{x}$
with vertex $\Psi(x)$ in $V$ such that $\Psi(U_x) $ is a neighborhood
of the vertex $\Psi(x) $ in $C_x$ and such that $\Psi|_{U_x}
:U_{x}\rightarrow C_{x}$ is an open map, where $C_x $ is endowed with
the subspace topology inherited from $V $. A map
$\Psi :X\rightarrow V$ is said to be {\bfi locally fiber connected}, if
for each
$x \in X $ there is an open neighborhood $U_x$ of $x$ such that
$\Psi^{-1}(\Psi(u))\cap U_{x}$ is connected for all $u\in
U_{x}$.  In Definitions~\ref{definition: local
convexity data} and~\ref{locally fiber connected condition (LFC)} we will elaborate
on these conditions. A continuous map between two topological spaces with Hausdorff
domain is said to be {\bfi proper} if it is closed and all its fibers are compact
subsets of its domain.

In spite of its generality this theorem can not be applied
to situations where the fibers $\Psi ^{-1}(v)$ are either not compact or
the map $\Psi $ is not closed because both conditions are
necessary for $\Psi $ to be a proper map. This is, for example,
one of the difficulties in the (direct) proof of the convexity
theorem due to Flaschka and Ratiu \cite{flaschkaratiu}
that contains as an important particular instance the case
of Poisson-Lie group actions on compact symplectic manifolds.
Alekseev \cite{alekseev} reproved the Poisson-Lie group action
result mentioned above by other means, reducing it to the case
covered by the Guillemin-Sternberg-Kirwan convexity result. His
method strongly uses the structure of Poisson-Lie groups and it
is not known how to extend it to other types of actions.

The main purpose of this paper is finding the essential features underlying all the
results mentioned above that ensure convexity. As we will see these properties
are of topological nature; more specifically, we will show that convexity is rooted on
the map being open onto its image and having local convexity data. These conditions
happen to be satisfied in all the classical convexity theorems that we discussed
previously and hence they can in principle be obtained as corollaries of the following
general theorem that can be found in Section~\ref{Openness and local convexity}:

\medskip

\noindent {\bfi Theorem}
\textit{ Let $f:X\rightarrow V$ be a continuous map from a
connected Hausdorff topological space $X $ to a Banach space $V$
that is open onto its image and has local convexity data.
Then the image $f(X)$ is locally convex. If, in addition,
$f(X)$ is closed in $V$ then it is convex.}

\medskip

Note first that $V$ is allowed to be infinite dimensional. Second, unlike the
local convexity data condition that can be found in~\cite{dazord, hilnebplank}, we do
not assume the cones $C_x $ to be
closed since it is not a reasonable assumption in infinite dimensions.

In the light of the result above, the convexity problem reduces to
giving necessary and sufficient conditions for a map that has
local convexity data to be open onto its image. In the paper (Section~\ref{Openness
and local convexity for momentum maps}) we will provide those characterizations for the
momentum maps associated to compact Lie group actions on symplectic manifolds. We will
split the problem in two cases:  when the map has connected fibers and when it has
only the locally fiber connectedness property. We shall also show that the openness of
the momentum map can be determined just by looking at its image and we will
illustrate this with two examples.

It is worth noting at this stage that Montaldi and Tokieda \cite{montaldi}
proved that the openness of the momentum map (relative to its
image endowed with the subspace topology) implies persistence of
extremal relative equilibria under every perturbation of the
value of the momentum map, provided the isotropy subgroup of
this value is compact. So the openness property of the momentum
map onto its image has interesting dynamical consequences.

Section~\ref{Openness and local convexity} contains a generalization of the
Lokal-global-Prinzip in~\cite{dazord, hilnebplank} that only requires the map to be
closed and to have a normal topological space as domain, instead of using the
properness condition. This degree of generality is needed to obtain convexity
directly from the Lokal-global-Prinzip in some of the examples that we present
and that generalize various results in the literature. For instance, in
Theorem~\ref{prato theorem better} we extend a  result of
Prato~\cite{prato} where we only require the properness of a single component of
the momentum map to conclude convexity.  
Additionally,  using our generalization of the 
Lokal-global-Prinzip we are able to drop in Section~\ref{Convexity for Poisson
actions of compact Lie groups}  the compactness hypothesis on the manifold in
the Flaschka-Ratiu convexity result~\cite{flaschkaratiu}.

It is worth mentioning that the generalization of the
Lokal-global-Prinzip presented in the paper includes infinite-dimensional situations.
This suggests that one could, in principle, use this tool in
dealing with convexity problems such as those in the papers of
Bloch, Flaschka, and Ratiu~\cite{bfr} or of
Neumann~\cite{neumann}. The implementation of this idea is not
free of difficulties and remains an open problem. This is due to
the lack of a Marle-Guillemin-Sternberg normal form in the
infinite dimensional setting which makes the local convexity data
property very difficult to check.

\section{Openness and local convexity}
\label{Openness and local convexity}

In this section we will study topological properties of maps
that have \textit{local convexity data}, a notion that we will
introduce below. This property holds typically for momentum
maps and, under certain supplementary topological conditions on
the map and on its domain of definition, it implies that the map
is open onto its image. We will also show that local convexity
of the image of a map that has local convexity data is implied
by openness of that map.

\paragraph{The Klee Theorem.} The passage from local convexity to
convexity is given by a classical result of Klee \cite{klee}
which we now present.

\begin{definition}
\label{definition: basic convexity}
Let $V $ be a topological vector space. If $x, y \in V $,
the {\bfi straight line segment}, or simply {\bfi segment}, $[x,
y]$ is defined by $[x,y]: = \{(1-\lambda)x + \lambda y \mid 0
\leq \lambda \leq 1 \}$. A subset $X \subset V $ is said to be
{\bfi convex\/}, if for any $x, y \in X $ we have $[x,y]\subset
X$. A subset $Y\subset V$ is called {\bfi locally convex\/} if
each point $y \in Y$ has a neighborhood $V_{y}$ whose intersection
with $Y$ is convex. A {\bfi polygonal path} is a continuous path
that is the union of segments. A subset $X \subset V$ is said to
be {\bfi polygonally connected\/} if any two points can be joined
by a polygonal path lying entirely in $X$.
\end{definition}

The following lemma is due to Kakutani and Tukey and will be
used in the proof of Klee's Theorem that insures the passage from
local to global convexity.

\begin{lemma}
\label{Kakutani}
In a topological vector space $V$, a connected locally convex set
$X$ is polygonally connected.
\end{lemma}

\noindent\textbf{Proof.\ \ }
Let $p$ be an arbitrary point of $X$ and let $X_{p}$ be the
set of all points of $X$ which can be joined to $p$ by a
polygonal path. The strategy is to prove that $X_{p}$ is both
an open and closed subset of $X$. The result then follows by
connectivity of $X$.

(i) $X_{p}$ is open in $X$: Take an arbitrary point $x\in X_{p}$
and a neighborhood $U_{x}$ of $x$ in $X$ chosen such that
$U_{x}\cap X$ is convex. This is possible since $X$ is locally
convex. Each $y \in U_x \cap X$ can be joined to $p $ by the
polygonal path obtained by adding to the polygonal path joining
$p$ to $x$ the straight line segment from $x$ to $y$ that is
guaranteed to lie entirely in $U_x \cap X$. Thus $p$ can be joined
to $y \in U_x$ by a polygonal path that lies entirely in $X$ which
proves that $y \in X_p $. Since $y \in U_x \cap X$ was arbitrary,
this shows that the open set $U_x \cap X $ in the relative
topology of $X $ lies in $X_p$, that is, that $X_p $ is open in
$X$.

(ii) $X_{p}$ is closed in $X$:  We will show that
$\overline{X}_{p}\cap X=X_{p}$.  To prove the non-trivial
inclusion $\overline{X}_{p}\cap X \subset X_{p}$, take an element
$y\in \overline{X}_{p}\cap X$ and recall that for any neighborhood
$U_{y}$ of $y$ in the topological vector space $V$ we have
$U_{y}\cap X_{p}\neq \varnothing$. We can chose, by hypothesis,
$U_{y}$ such that $U_{y}\cap X$ is convex. Thus there is an
element $z\in U_{y}\cap X_p$ such that the straight line segment
joining $z $ to $y $ lies entirely in $U_y \cap X $. However, $z
\in X_p $, so there is a polygonal path that lies entirely in $X $
joining $p $ to $z $. Adding to this path the segment joining $z $
to $y $ yields a polygonal path lying entirely in $X $ that joins
$p $ to $y $, which proves that $y \in X_p $. \quad $\blacksquare$

\medskip

The next theorem is due to Klee \cite{klee} and gives the
connection between local convexity and convexity. Since this
result does not seem to be widely known, we shall reproduce its
proof.

\begin{theorem}[Klee]
\label{klee}
Each closed connected and locally convex subset of a
topological vector space is convex.
\end{theorem}

\noindent\textbf{Proof.\ \ }
For $x, y, z \in V$ denote by $\mathcal{CH}(x,y,z) : = \{ v \in V
\mid v = \lambda_1 x + \lambda_2 y + \lambda_3 z, \lambda_1,
\lambda_2, \lambda_3 \geq 0, \lambda_1 + \lambda_2 + \lambda_3 =
1 \}$ the convex hull of the set $\{x, y, z \}$.

Let $X$ be a closed connected and locally convex subset of
a topological vector space $V$. The strategy is to prove that we
can cut corners in $X$, i.e. if $[x,y] \subset X$ and
$[y,z]\subset X$ then $[x,z]\subset X$. To prove this, we first
show that if $[x,y]\subset X$ and $[y,z]\subset X$ with $x\neq
y$ then $\mathcal{C}\mathcal{H}\{q,y,z\}\subset X$ for some
$q\in \lbrack x,y)$ (that is, $q \in [x,y]$ and $q  \neq y$). Let
$K$ be the collection of all points $p\in \lbrack y,z]$ such
that there exists a point $q\in \lbrack x,y)$ depending on $p$
with $\mathcal{CH}\{q,y,p\}\subset X$. The set $K$ is non-empty
since $y\in K.$ Again we will prove that $K$ is open and closed
in $[y,z].$

\begin{figure}[htb]
\begin{center}
\includegraphics{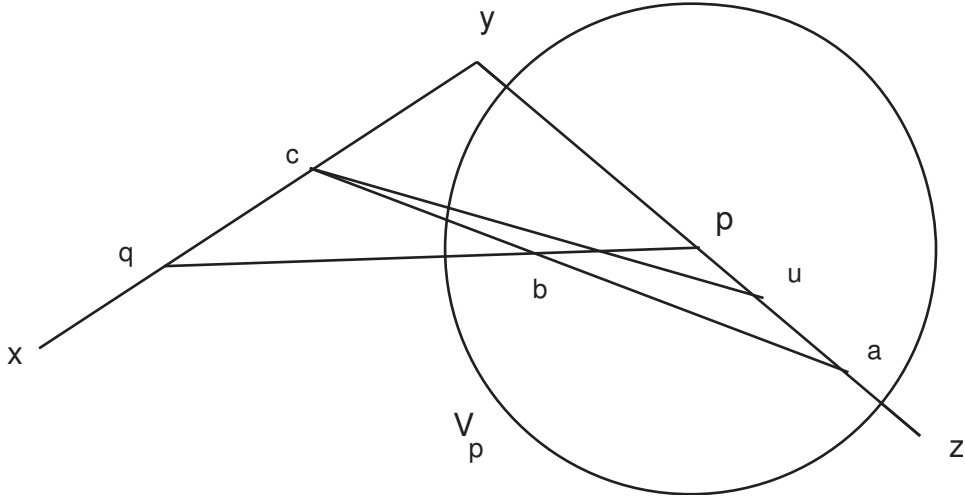}

\caption{The proof of openness of $K$}
\label{fig:orbsta}
\end{center}
\end{figure}

To prove that $K$ is open let $p \in K$ and $U_{p}$ an open
neighborhood of $p$ in $V$ such $V_{p}:=U_{p}\cap X$ is convex.
We have to consider two cases. First, assume that  $a\in
[ p,z]\cap V_{p}$ is arbitrary. Then  we can find an element
$b\in (p,q] \cap V_p$ such that
$\mathcal{CH}\{b,p,a\}\subset V_{p}$; this is possible by
convexity of $V_p $. The line that goes through the points $a$
and $b$ intersects $[x,y)$ in a point $c$. We have
$\mathcal{CH}\{c,y,a\}\subset X$; indeed,
$\mathcal{CH}\{c,y,a\}$ is the union of the tetragon whose
vertices are $c, y, p, b $ which is included in
$\mathcal{CH}\{q,y,p\}$ by construction and is hence a subset of
$X $ and the triangle $\mathcal{CH}\{b, p, a\} \subset V_p
\subset X$ as was just shown. Since $\mathcal{CH}\{c, y, a\}
\subset X$, it follows that $\mathcal{CH}\{c,y,u\} \subset
\mathcal{CH}\{c,y,a\} \subset X$ whenever
$u\in [y,a]$. In particular, this shows that for any $u \in U_p
\cap [y,z] = V_p \cap [y,z]$ we have $\mathcal{CH}\{c, y, u\}
\subset X$, that is, $u \in K$. In particular, the arbitrarily
chosen $a \in [p,z] \cap V_p $ is in $K $. The second case is if
$a \in [y,p] \cap V_p$ is arbitrary. Then $\mathcal{CH}\{q, y,
a\} \subset
\mathcal{CH}\{q, y, p\} \subset X $ which shows that $a \in K $.
So for all $p \in K $ we can find an open interval around $p $ on
$[y,z]$ which is included in $K $, which shows that $K $ is open
in $[y,z]$.

To prove that $K$ is closed in $[y,z]$ it will be shown that
$\overline{K}\cap [y,z] = K\cap [y,z]$. To prove the non-trivial
inclusion, let $p\in \overline{K}\cap [y,z]$ and $U_{p}$ an open
neighborhood of $p$ in $V$ such that $V_{p}:=U_{p}\cap X$ is
convex. Then $U_p \cap K = V_p \cap K \neq \varnothing $, so let
$a \in V_p \cap K$. Thus, there is some $q  \in [x, y)$ such
that $\mathcal{CH}\{q,y,a\} \subset X$. Note that if $a \in [p,
z]$ then  $p \in [y, a]$. Since $a \in K$ this implies that $p
\in K
$. So assume that $a \in [y,p]$. Since $a \in K $, there is some
$q \in [x, y)$ such that $\mathcal{CH}\{q, y, a\} \subset X $. By
convexity of
$V_p$, for any $b\in (a,q] \cap V_p$ we have
$\mathcal{CH}\{b,a,p\}\subset V_{p}.$ But then the line that
goes trough the points $p$ and
$b$ intersects $[x,y)$ in a point $c$. We have
$\mathcal{CH}\{c,y,p\}\subset \mathcal{CH}\{q,y,a\}\cup
\mathcal{CH}\{b,a,p\} \subset X$ which implies that $p\in K$.
This shows that in either case $p \in K $ which in turn proves
that $K$ is closed in $[y,z]$.

\begin{figure}[htb]
\begin{center}
\includegraphics{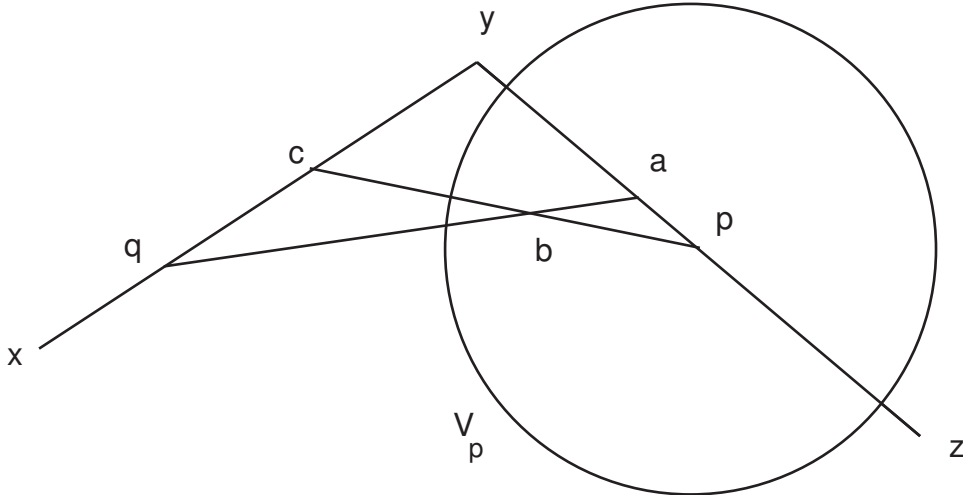}
\caption{The proof of closedness of $K$}
\label{fig:orbsta}
\end{center}
\end{figure}

Since $K$ is both open and closed in $[y,z]$ we have $K=[y,z]$.
Thus,
\begin{equation}
\label{relation just proved}
\text{if } [x,y]\subset X \text{ and } [y,z]\subset X \text{ with
} x\neq y \text{ then } \mathcal{CH}\{q,y,z\}\subset X \text{ for
some } q\in [x,y).
\end{equation}
Now
we will prove that if
\begin{equation}
\label{klee relation}
[x,y]\subset X \quad \text{and} \quad [y,z]\subset X \quad
\text{then} \quad [x,z]\subset X.
\end{equation}
To see this, assume that $[x,y]\subset X$  and $[y,z]\subset X$
and define $M$ to be the collection of all points $p\in [x,y]$
such that $\mathcal{CH}\{p,y,z\}\subset X$. The set $M$ is
nonempty since $y\in M$ and it is closed because $X$ is closed in
$V$. Indeed, let $p_n \in M $ be a sequence in $M$ that converges
to a point $p \in [x,y]$. We need to show that $p \in M $. Since
$p_n \in M $ we have that $\mathcal{CH}\{p_n,y, z\} \subset X $
and hence $\cup_{n \in \mathbb{N}} \mathcal{CH}\{p_n,y,z\} \subset
X $. Therefore $p \in \overline{ \cup_{n \in \mathbb{N}}
\mathcal{CH}\{p_n,y,z\}} \subset \overline{X} = X $ since $X $ is
closed.

If we show that $M$ is also open, then $M = [x,y]$ and so
\eqref{klee relation} is proved. Let
$p\in M$ with $p\neq x$. Then $[x,p]\subset [x, y] \subset X$ and
$[p,z]\subset \mathcal{CH}\{p,y,z\} \subset X$. By
\eqref{relation just proved}, there is a point $q\in [x,p)$ such
that $\mathcal{CH}\{q,p,z\}\subset X$. But then
$\mathcal{CH}\{a,y,z\}\subset X$ for each $a\in [y,p]$. Indeed,
if $a \in [p,y]$ then $\mathcal{CH}\{a,y,z \} \subset
\mathcal{CH}\{p,y,z \}  \subset X$. Thus, $[y,p] \subset M $. If
$a \in [q,p]$, then
$\mathcal{CH}\{a,y,z \} = \mathcal{CH}\{a,p,z \} \cup
\mathcal{CH}\{p,y,z\}  \subset \mathcal{CH}\{q,p,z \} \cup
\mathcal{CH}\{p,y,z\}  \subset X$. Thus, for every $a \in
[q,p]$, we have $\mathcal{CH}\{a,y,z\} \subset X $  which
means that $a \in M $. This proves that $(q, y] \subset M$.
Since  $(q, y]$ is an open neighborhood of $p $ in the relative
topology of $[x,y]$ this argument shows that $M $ is open in
$[x,y]$.

The convexity of $X$ is proved as follows. Consider two
arbitrary points $p$ and $q$ in $X$. By Lemma \ref{Kakutani} it
follows that there are points $x_{0}=p, x_{1}, \dots , x_{n}=q$
in $X$ such that $[x_{i-1},x_{i}]\subset X$ for every $i\in
\{1,...,n\}$. Applying $n-1$ times \eqref{klee relation} we
obtain $[p,q]\subset X$ which proves that $X$ is convex.
\quad $\blacksquare$

\begin{remark}
\label{Klee remark}
{\rm By an analogous proof we can replace the condition that
the set $X $ is closed in $V$ with the weaker condition that the
set is closed in a convex subset of $V$ endowed with the induced
topology from $V$. This remark will be used in Theorem
\ref{conv}.}
\end{remark}

\paragraph{Maps with local convexity data.} Let $f:X\rightarrow V$
be a continuous map defined on a connected Hausdorff topological
space
$X $ with values in a locally convex vector space $V$. On the
topological space $X$ define the following equivalence relation:
declare two points
$x,y\in X$ to be equivalent if and only if $f(x)=f(y)=v$ and
they belong to the same connected component of $f^{-1}(v)$.
The topological quotient space will be denoted by $X_{f}:=X/R$,
the projection map by $\pi _{f}:X\rightarrow X_{f}$, and the
induced map on $X_{f}$ by $\widetilde{f}:X_{f}\rightarrow V$.
Thus, $\widetilde{f} \circ \pi_f = f$ uniquely
characterizes $\widetilde{f}$. The map
$\widetilde{f}$ is continuous and if the fibers of $f$
are connected then it is also injective.

\medskip

The following elementary topological fact will be used
several times later on.

\begin{lemma}
\label{conex}
Let $f:X\rightarrow Y$ be a continuous map between two
topological spaces. Assume that $f $ has connected fibers and is
open or closed. Then for every connected subset $C$ of $Y$ the
inverse image $f^{-1}(C)$ is connected.
\end{lemma}

\noindent\textbf{Proof.\ \ } Suppose that $f $ is an open map,
$C \subset Y $ is connected, and $f ^{-1} (C)$ is not
connected. Then there exist two open sets $U_1, U_2$ in
$X$ such that $f ^{-1} (C) = (U_1 \cap f ^{-1} (C)) \cup (U_2
\cap f ^{-1} (C))$, $ U_1 \cap f ^{-1} (C) \neq \varnothing $,
$U_2 \cap f ^{-1} (C) \neq \varnothing $,  and $U_1 \cap U_2 \cap
f ^{-1} (C) = \varnothing $. Note that $C = f(f ^{-1} (C)) = f(
(U_1\cap f ^{-1} (C)) \cup (U_2 \cap f ^{-1} (C))) = f (U_1
\cap f ^{-1} (C)) \cup f (U_2 \cap f ^{-1} (C)) \subset
(f(U_1) \cap C)\cup (f(U_2) \cap C)$. Conversely, since $f(U_1)
\cap C \subset C $ and $f(U_2) \cap C \subset C $ it follows
that $(f(U_1) \cap C) \cup (f(U_2) \cap C) \subset C $ which
proves that $(f(U_1) \cap C) \cup (f(U_2) \cap C) = C$. Also,
$f(U_1) \cap C \supset f(U_1 \cap f ^{-1} (C)) \neq \varnothing$
and $f(U_2) \cap C \supset f(U_2 \cap f ^{-1} (C)) \neq
\varnothing$. By openness of $f$, the sets $f(U_1)$ and
$f(U_2)$ are open in $Y$ so that connectedness of $C $ implies
that $f(U_1) \cap f(U_2) \cap C \neq \varnothing $.

If $c \in f(U_1) \cap f(U_2) \cap C$ then $f ^{-1}(c) = (U_1 \cap
f ^{-1}(c)) \cup (U_2 \cap f ^{-1}(c))$. The inclusion
$\supset $ is obvious. To prove the reverse inclusion $\subset$,
let $x \in f ^{-1}(c) \subset f ^{-1}(C)$. Thus $x \in U_1 \cap
f ^{-1}(C) $ or $x \in U_2 \cap f ^{-1}(C)$. Since $x \in
f^{-1}(c) $ by hypothesis, this implies that $x \in U_1 \cap
f^{-1}(c) $ or $x \in U_2 \cap f ^{-1}(c)$ which proves the
inclusion $\subset $. Note also that $U_1 \cap f ^{-1}(c) \neq
\varnothing $ since $c \in f (U_1)$. Similarly, $U_2 \cap f
^{-1}(c) \neq \varnothing $. Finally, $U_1 \cap U_2\cap f
^{-1}(c)  \subset U_1 \cap U_2 \cap f ^{-1}(C) = \varnothing $.
Thus the fiber $f ^{-1}(c)$ can be written as the disjoint
union of the two open nonempty sets $U_1 \cap f ^{-1}(c)$ and
$U_2 \cap f ^{-1}(c)$, which contradicts the connectedness
hypothesis of the fibers of $f$.

The proof for $f$ a closed map is identical to the one above by
repeating the same argument for $U_1 $ and $U_2 $ closed subsets
of $X$.
\quad $\blacksquare$

\begin{definition}
\label{definition of the cone} Let $V$ be a topological vector
space. A set $C \subset V $ is called a {\bfi cone} with {\bfi
vertex\/} $v_0$ if for each $\lambda \geq 0$ and for each $v \in
C$, $v \neq v_{0}$, we have $(1-\lambda)v_{0}+\lambda v \in C$. If
the set $C $ is, in addition, convex then $C $ is called a {\bfi
convex cone}. Note that, by definition, the vertex $v_0 \in C$.
\end{definition}

\begin{definition}
\label{definition: local convexity data} 
The continuous map
$f:X\rightarrow V$ defined on a connected locally connected
Hausdorff topological space $X$ with values in a locally convex
topological vector space $V$ is said to have \textbf {local
convexity data} if for each $x\in X$ and every sufficiently small
neighborhood $U_{x}$ of $x$ there exists a convex cone
$C_{x,f(x),U_{x}}$ in $V$ with vertex at $f(x)$ such that
\begin{description}
\item[(VN)] $f(U_{x})\subset C_{x,f(x),U_{x}}$ is a neighborhood
of the vertex $f(x)$ in the cone $C_{x,f(x),U_{x}}$ and

\item[(SLO)]$f|_{U_{x}}:U_{x}\rightarrow C_{x,f(x),U_{x}}$ is an
open map and for any neighborhood $U_{x}^{^{\prime
}}$ of $x$, $U_{x}^{^{\prime
}}\subset U_{x}$, the set $f(U_{x}^{^{\prime }})$ is a neighborhood
of the vertex $f(x)$ in the cone $C_{x,f(x),U_{x}}$,
\end{description}
where the cone $C_{x,f(x),U_{x}}$ is endowed with
the subspace topology induced from $V$.
If the associated cones  $C_{x,f(x),U_{x}}$ are such that
$C_{x,f(x),U_{x}} \cap f(X)$ is closed in $f(X)$, then we say that
$f$ has {\bfi local convexity data with closed cones}.
\end{definition}

\begin{remark}
\normalfont 
We are using the (VN) condition as an abbreviation for
``vertex neighborhood condition" and the (SLO) condition as an
abbreviation for ``strong local openness condition".
Note that in the
case when the associated cones $C_{x,f(x),U_{x}}$ are closed in
$f(X)$ the second condition in (SLO) is automatically implied by
the openness of $f|_{U_{x}}:U_{x}\rightarrow C_{x,f(x),U_{x}}$.
\end{remark}

\begin{remark}
\label{first LC remark} \normalfont Let $V$ be a locally convex
topological vector space, $C_{1}$ and $C_{2}$ two cones in $V$
with vertex at zero, and $V_0$ a neighborhood of zero in $V$.
Suppose that $C_{1} \cap C_2 \cap V_{0}$ is a neighborhood of zero
in $C_2$. Then $C_{2}\subset C_{1}$.

Indeed, since scalar multiplication is a continuous operation
with respect to the subspace topology on $C_{2}$ induced by
the topology of $V$ and since $V $ is a locally convex
topological space, we obtain that for every
$x\in C_{2}$ there exists some $\lambda>0$ such that $\lambda x
\in C_{1} \cap C_2 \cap V_{0}\subset C_1$. Since $C_1$ is a
cone, it follows that
$t\lambda x \in C_1$ for every $t>0$ and hence $x \in C_1$.

By translation, the same property holds if the common vertices of
the cones $C_1$ and $C_2$ are at some other point of $V$. We shall
use this observation several times in the remarks that follow.
\end{remark}

\begin{remark}
\label{second LC remark} \normalfont $C_{x,f(x),U_{x}}$ does not
depend on the neighborhood $U_{x}$ in the sense that if
$U_{x}^{^{\prime }}\subset U_{x}$ is another neighborhood of $x$
then $C_{x,f(x),U_{x}^{^{\prime }}}=C_{x,f(x),U_{x}}$.

Indeed, we have $f(U_{x}^{^{\prime}})\subset
C_{x,f(x),U_{x}^{^{\prime}}}$ is a neighborhood of the vertex
$f(x)$ in $C_{x,f(x),U_{x}^{^{\prime}}}$, $f(U_{x})\subset
C_{x,f(x),U_{x}}$ is a neighborhood of the vertex $f(x)$ in
$C_{x,f(x),U_{x}}$, and $f(U_{x}^{^{\prime
}})=V^{^{\prime}}_{f(x)}\cap C_{x,f(x),U_{x}^{^{\prime}}}\subset
f(U_{x})=V_{f(x)}\cap C_{x,f(x),U_{x}}$, where
$V_{f(x)}^{^{\prime}}$ and $V_{f(x)}$ are two open neighborhoods
of $f(x)$ in $V$. By the argument used in Remark \ref{first LC
remark} it follows that $C_{x,f(x),U_{x}^{^{\prime }}}\subset
C_{x,f(x),U_{x}}$. The (SLO) condition shows that
$f(U_{x}^{^{\prime }}) = f(U_{x}^{^{\prime }}) \cap
C_{x,f(x),U_{x}^{^{\prime}}} \cap C_{x,f(x),U_{x}}$ is a
neighborhood of the vertex $f(x)$ in $C_{x,f(x),U_{x}}$. Again
Remark \ref{first LC remark} implies that $C_{x,f(x),U_{x}}
\subset C_{x,f(x),U_{x}^{^{\prime}}}$.

Thus, since $C_{x,f(x),U_{x}^{^{\prime}}}$ is independent of the
neighborhood $U_x$, we shall write $C_{x,f(x)}$ for the cone in
Definition \ref{definition: local convexity data}.
\end{remark}

\begin{remark}
\label{third LC remark} \normalfont In Remark \ref{fourth LC
remark} we shall need the following statement (see e.g. \cite{ChVo} or
\cite{HoYo}): a topological space is connected if
and only if every open covering $\{U_{\alpha }\mid \alpha \in
\mathcal{A}\}$ has the property that for each pair of points $x$,
$y$ there exists a finite sequence $\{\alpha_{1}$, ..., $\alpha
_{k} \} \subset \mathcal{A}$ such that $x\in U_{\alpha _{1}}$,
$y\in U_{\alpha _{k}}$, and $U_{\alpha_{i}}\cap
U_{\alpha_{i+1}}\neq \varnothing$ for all $i=1, \dots , k-1$. This
finite family of open sets $\{U_{\alpha _{1}}, \dots, U_{\alpha
_{k}} \}$ is called a {\bfi finite chain\/} linking $x$ to $y$. A
topological space that has the property that for each open cover
and arbitrary points $x,y $ there is a finite chain formed by
elements of this open cover linking $x $ to $y $ is also called
{\bfi well-chained\/}. The statement above asserts hence that a
topological space is connected if and only if it is well-chained.

\end{remark}

\begin{remark}
\label{fourth LC remark}
\normalfont
$C_{x,f(x)}$ depends only on the connected components of
$f^{-1}(f(x))$, that is, if $y$ is in the same connected
component of $f^{-1}(f(x))$, then $C_{x, f(x)} = C_{y, f(y)}$.

Let $U_x $ be the open neighborhood of $x \in X $ guaranteed by
Definition \ref{definition: local convexity data}. Let $y\in
f^{-1}(f(x))$ be in the same connected component of $f^{-1}(f(x))$
as $x$ and $y \in U_x$. For this $y $ let $U_y $ be the open
neighborhood of $y$ in  Definition \ref{definition: local
convexity data}  which we can choose such that $U_{y}\subset
U_{x}$. The fact that $f(U_{y})\subset f(U_{x})$ and $f(U_{y})$ is
a neighborhood of $f(x)$ in $C_{y,f(x)}$ and $f(U_{x})$ is a
neighborhood of $f(x)$ in $C_{x,f(x)}$ implies that
$C_{y,f(x)}\subset C_{x,f(x)}$. 
For the reverse inclusion observe
that $f(U_{y})\subset C_{y,f(x)} $ is a neighborhood of $f(x)$ and
by the (SLO) condition we also have that $f(U_{y})$ is open in
$C_{x,f(x)}$ which shows that $f(U_y) \cap  C_{x,f(x)} = f(U_y)
\cap C_{x,f(x)} \cap C_{y,f(x)} $ is an open neighborhood of
$f(x)$ in $C_{x,f(x)}$.
By Remark \ref{first LC remark}, we obtain
$C_{x,f(x)} \subset C_{y,f(x)}$. Thus for any $y \in U_x$ that
also lies in the same connected component of $f ^{-1}(f(x))$ as $x
$, we have $C_{x,f(x)} = C_{y,f(x)}$.

So, on the connected component $E_x$ of the fiber $f
^{-1}(f(x))$ containing $x $, using Remark \ref{third LC
remark}, we obtain $C_{x,f(x)} = C_{y,f(x)}$ for any $y \in
E_x$. Thus, if the fibers $f^{-1}(f(x))$ are connected
we can erase $x$ from $C_{x,f(x)}$.
\end{remark}

Our strategy to prove local convexity for the image of a map that
has local convexity data is to prove that it is open onto its
image.

\begin{theorem}
\label{conv}
Let $X $ be a connected Hausdorff topological space, $V $ a
locally convex topological vector space, and $f:X\rightarrow V$
a continuous map that has local convexity data. If
$f$ is open onto its image then $f(X)$ is a locally convex subset
of $V$. Moreover, if $f(X)$ is closed in a
convex subset of $V$ then it is convex.
\end{theorem}

\noindent\textbf{Proof.\ \ } Let $v\in f(X)$ be arbitrary and take
$x\in f^{-1}(v)$. By the condition (VN), there exists a
neighborhood $U_{x} \subset X$ of $x$ such that $f(U_{x}) \subset
C_{x,f(x)}$ is open in $f(X)$; $C_{x,f(x)}$ is the convex  cone
with vertex at $v=f(x)$ given in Definition \ref{definition: local
convexity data}. Thus, shrinking $U_x$ if necessary, using
condition (SLO), and the local convexity of the topological vector
space $V $, we can find a convex neighborhood $V_v $ of $v$ in $V$
such that $f(U_{x})=V_{v}\cap C_{x,f(x)}$. Since $f $ is open onto
its image, the neighborhood $V_v $ can be shrunk further to a
convex neighborhood of $v$, also denoted by $V_v$, such that $f
(U_x) = V_v \cap f(X)$. Taking this as the neighborhood of $v $
and shrinking $U_x$ if necessary, we get $V_{v}\cap C_{x,f(x)} =
f(U_x) =  V_v \cap f(X)$. Since the intersection of two convex
sets is convex, it follows that $V_v \cap C_{x,f(x)}$ is also
convex. Thus the point $v \in f(X)$ has a neighborhood $V_v $ in
$V $ such that $V_v \cap f(X) $ is convex, that is, $f(X)$ is
locally convex (see Definition \ref{definition: basic convexity}).

Now assume, in addition, that $f(X)$ is closed in a convex subset $C$ of $V$.
Connectedness of $X$ and  continuity of $f $ imply that $f(X)$ is connected.
Therefore, $f(X)$ is a closed, connected, and locally  convex subset of $C$ (by what
was just proved), so Klee's theorem (see Theorem~\ref{klee}) ensures that $f(X)
$ is convex in $C$. Since $C$ is convex in $V$, this implies that $f(X)$ is convex in $V$.
\quad $\blacksquare$

\medskip

\begin{remark}
\label{remark for open maps}
\normalfont
Note that in Theorem \ref{conv} it was not assumed that the cones
that give the local convexity data are closed. This could be
useful for convexity theorems with infinite dimensional range.
It is also worth noting that connectivity of $X $ was not really
used. In other words, the theorem above holds for each connected
component of $X $.
\end{remark}

It is not obvious when a map that has local convexity data is open
onto its image. In order to give sufficient conditions under which
this happens we need a few more preliminary results. First we need
the following concept.

\begin{definition}
\label{locally fiber connected condition (LFC)}
Let $X $  and $Y$  be two topological spaces and $f:X
\rightarrow  Y$  a continuous map. The subset $A \subset X$
satisfies the {\bfi  locally fiber connected condition (LFC)}
if $A$ does not intersect two different connected components of
the fiber
$f ^{-1}(f (x))$, for any $x \in A  $.

Let $X $ be a connected, locally connected, Hausdorff topological
space and $V $ a locally convex topological vector space.  
The
continuous map $f: X \rightarrow V$ is said to be \textbf{locally
fiber connected} if for  each $x\in X$, any open neighborhood of
$x $ contains a neighborhood $U_{x}$ of $x$ such that
$U_{x}$ satisfies the {\bf (LFC)} condition.
\end{definition}

Recall that $X_f$ denotes the quotient topological space of $X$
whose points are the connected components of the fibers of the
continuous map $f:X \rightarrow V$ and $\widetilde{f}:X _f
\rightarrow  V $  is the map such that $\widetilde{f} \circ \pi
_f = f $.
Note that a subset  $A \subset X$
satisfies  (LFC) if and only if $\widetilde{f}|_{\pi _f (A)} $
is injective. Similarly, $f$ is locally fiber connected if and
only if for any $x \in  X $, any open neighborhood  of $x
$  contains an open neighborhood $U_x  $  of $x$ such that the
restriction of $\widetilde{f} $  to $\pi _f(U _x)$  is injective.

\begin{lemma}
\label{primatop} Assume that  the continuous map $f:X\rightarrow
V$ has local convexity data and is also locally fiber connected.
\begin{itemize}
\item[{\rm \textbf{(i)}}] The quotient projection $\pi_{f}: X
\rightarrow X_f$ is an open map.

\item[{\rm \textbf{(ii)}}] If $f$ has connected fibers then
$X_{f}$ is Hausdorff.
\end{itemize}
\end{lemma}

\noindent\textbf{Proof.\ \ } \textbf{(i)} We begin by noting that 
it suffices to prove that for each $x \in
X$ and for each open connected neighborhood $U_x$ that satisfies
(VN), (SLO), and (LFC), its image $\pi_f(U_x)$ is open in $X_f$.
This is so because any open neighborhood of $x $ contained in $U_x
$ also satisfies the same three  properties.

Take $x$ an arbitrary point in $X$ and $U_x$ a connected
neighborhood of $x$ that satisfies (VN), (SLO), and (LFC). We will
prove that $\pi_f^{-1}(\pi_f(U_x))$ is open in $X$. To show this,
let $y$ be arbitrary in $\pi_f^{-1}(\pi_f(U_x))$. Then the
connected component $E_y$ of $f ^{-1}(f(y)) $ that contains $y $
intersects $U_x$. Let $y' \in E_y \cap U_x$. Choose a neighborhood
$U_{y'}\subset U_x$ that also satisfies (VN), (SLO) and (LFC).
Since $E_y$ is a connected space with respect to the induced topology from $X$, 
it is well chained with respect to this
topology, hence we can find a finite chain of open
sets in $X$ satisfying (VN), (SLO) and (LFC) such that
$U_1=U_{y'}$ as above, $y\in U_n$ and $U_i\cap U_{i+1}\cap
E_y\neq\varnothing$.  The set
$W_{yy'}:=\bigcap_{i=1}^{n-1} f(U_i\cap U_{i+1})$ is an open
subset of $C_{y,f(y)}$ because of (SLO) and due to the fact that the
associated cones depend only on the connected components of the
fibers. Additionally, $W_{yy'}$ is nonempty since by construction it contains
the point $f(E_y)$, because  $U_i\cap
U_{i+1}\cap E_y\neq\varnothing$, for any $i \in \{1, \ldots, n-1\}$. We now show
that $O_y:=U_n\cap f^{-1}(W_{yy'})$ is an open subset of $\pi_f^{-1}(\pi_f(U_x))$
that contains $y$. Indeed, $O_y$ is open in $X$ because it is an open subset of
$U_n$.  By the construction of $O _y$, $y$ clearly belongs to $O_y$. It remains to be
shown that $O _y \subset  \pi _f^{-1} (\pi _f(U _x))$.

Let $z\in O_y$
arbitrary. The connected component $E_z$ of $f^{-1}(f(z))$ that
contains $z$ intersects $U_n$. By the construction of $O_y$, the
fiber $f^{-1}(f(z))$ intersects every open set $U_i\cap
U_{i+1}$ of the finite chain that links $y$ and $y'$. 
The (LFC) property guarantees that the connected component of $f^{-1}(f(z))$ that
intersects $U_n\cap U_{n-1}$ has to be $E_z$. Repeating the
argument we obtain that $E_z\cap U_{y'}\neq \emptyset$ and
$E_z\cap U_{y'}\subset U_{x}$. Then $E_{z}\subset
\pi_f^{-1}(\pi_f(U_x))$ for every $z\in O_{y}$, which proves that
$O_y\subset\pi_f^{-1}(\pi_f(U_x))$ and, consequently,
$\pi_f^{-1}(\pi_f(U_x))$ is an open subset of $X$.

\noindent\textbf{(ii)} Since all the fibers of $f $ are
connected and $f $ is continuous, it follows that the graph of the
relation determined by $f$ is closed. As $\pi_{f}$ is open we
obtain that $X_{f}$ is a Hausdorff space.
\quad$\blacksquare$

\medskip

\begin{remark}
\normalfont
Notice that neither local convexity data nor the locally fiber
connected condition alone would imply openness of $\pi_f$. Indeed,
if we consider the example where $X$ is the square in ${\mathbb
R}^{2}$ with vertices $(2,2)$, $(-2,2)$, $(-2,-2)$, $(2,-2)$ minus
the interior of the square with vertices $(0,1)$, $(-1,0)$,
$(0,-1)$, $(1,0)$, and $f(x,y)=y$, then $\pi_f$ is not open.
However, $f$ has local convexity data but it does not satisfy the
locally fiber connected condition. If we consider the example
where $X$ is the square in ${\mathbb R}^{2}$ with vertices
$(2,2)$, $(-2,2)$, $(-2,-2)$, $(2,-2)$ minus the interior of the
rotated square with vertices $(1,1)$, $(-1,1)$, $(-1,-1)$,
$(1,-1)$, and $f(x,y)=y$, then again $\pi_f$ is not open and $f$
satisfies the locally fiber connected condition but does not have
local convexity data, precisely because it does not satisfy the
(SLO) condition.
\end{remark}

An immediate consequence of Lemma \ref{primatop} and Theorem
\ref{conv} is the following corollary. 
It states that for a map
that has local convexity  data and connected fibers the
condition to be open onto its image, and consequently to have a
locally convex image,  is implied by the condition to be closed
onto its image or by the stronger condition to be a proper map.

\begin{corollary}
\label{imediat}
Let $X $ be a connected, locally connected, Hausdorff topological space, $V $ a locally convex 
topological vector space, and $f : X \rightarrow V$ a continuous map that has local convexity data. 
Assume that $f $ is a closed map onto its image $f(X)$ and that it has connected fibers. Then $f $ is 
open onto its image $f(X)$ and $f(X)$ is locally convex. Moreover, if $f(X)$ is closed then it is
convex.
\end{corollary}

\noindent\textbf{Proof.\ \ }
The hypothesis implies that the  induced map
$\widetilde{f}:X_{f}\rightarrow {f(X)}$, uniquely determined by
the equality $\widetilde{f} \circ \pi_f = f$, is a homeomorphism.
Indeed, closedness of $\widetilde{f} $ follows from the identity
$\widetilde{f}(A) = f (\pi_f^{-1}(A))$ for any subset $A$ of
$X_f$. 
Since $\widetilde{f} $ is open onto $f(X)$ it follows
that $f = \widetilde{f} \circ \pi_f $ is also open onto its
image. The rest is a consequence of Theorem \ref{conv}.
\quad $\blacksquare$

\medskip

Note that by Lemma \ref{conex}, a necessary condition for the
map $f$ that has connected fibers to be open onto its image is
that the inverse image of any connected set in $f(X)$ is
connected in $X$. The next proposition states that if $f $ has
local convexity data with closed cones, this condition is also
sufficient.

\begin{proposition}
\label{criteriu}
Let $f:X\rightarrow V$ be a continuous map that has local
convexity data with closed cones. If the fibers of $f$ are
connected and for every point $v\in f(X)$ and for all small
neighborhoods $V_{v}$ of $v$ the set $f^{-1}(V_{v})$ is
connected, then $f$ is open onto its image.
\end{proposition}

\noindent\textbf{Proof.\ \ } Suppose that $f$ is not open onto its
image, has local convexity data with closed cones, and connected
fibers. So there exists a point $x\in X$ and an open neighborhood
$U_x $ (included in a neighborhood of $x $ from the definition of
local convexity data) such that $V_{f(x)} \cap C_{x,f(x)} = f(U_x)
\subsetneq f(X) \cap V_{f(x)}$ for some open neighborhood
$V_{f(x)}$ of $f(x) $ in $V$. Consequently, $(V_{f(x)}\cap
f(X))\backslash f(U_{x}) \neq \varnothing$ is open in $f(X)$ since
$f(U_x) = V_{f(x)} \cap C_{x,f(x)} \cap f(X) $ is closed in the
topology of $V_{f(x)} \cap f(X) $ due to the  fact that
$C_{x,f(x)} \cap f(X)$ is closed in $f(X)$ (by the closed cone
hypothesis). We can also choose
$V_{f(x)}$ small enough so that the hypothesis holds for it, that
is, $f^{-1}(V_{f(x)}\cap f(X)) $ is connected.

Note that connectedness of the fibers, and thus bijectivity
of $\widetilde{f}$, implies that
$\widetilde{f}^{-1}(f(A)))=\pi_{f}(A)$ for any subset $A$ of $X$.

The sets that enter in the equality
\begin{equation*}
\widetilde{f}^{-1}(V_{f(x)}\cap f(X))
= \widetilde{f}^{-1}((V_{f(x)}\cap f(X))\backslash
f(U_{x}))\cup \widetilde{f}^{-1}(f(U_{x}))
\end{equation*}
or equivalently
\begin{equation*}
\pi _{f}(f^{-1}(V_{f(x)}\cap f(X)))=\widetilde{f}^{-1}((V_{f(x)}
\cap f(X))\backslash f(U_{x}))\cup \pi _{f}(U_{x})
\end{equation*}
are all open because $\pi _{f}$ is open and we also have that
$\widetilde{f}^{-1}((V_{f(x)}\cap f(X))\backslash
f(U_{x}))\cap \widetilde{f}^{-1} (f(U_{x}))=\varnothing$.
But this contradicts the connectivity of $\pi _{f}
(f^{-1}(V_{f(x)}\cap f(X)))$.
\quad $\blacksquare$

\medskip

\paragraph{The Lokal-global-Prinzip.}The rest of this section is
dedicated to the generalization of the ``Lokal-global-Prinzip''
used in the convexity proof of Condevaux, Dazord and Molino
\cite{dazord} and of Hilgert, Neeb, and Plank
\cite{hilnebplank}. On one hand we will prove that one can drop
the condition on the compactness of the fibers necessary in the
classical proof and still have the same conclusions and on the
other hand we will extend this result to maps that have an
infinite dimensional target. In order to avoid  the compactness
condition (we will call a topological space compact if it is Hausdorff and every
open cover has a finite subcover) on the fibers we need to further investigate the
topology of the problem.  We start with a short account of some topological results
needed in the sequel.

\begin{theorem}[\cite{engelking}]
\label{caracterizare inchise}
Let $f:X \rightarrow Y $ be a continuous mapping.
\begin{itemize}
\item[{\rm \textbf{(i)}}] $f$ is closed
if and only if for every $B\subset Y$ and every open set $A\subset
X$ which contains $f^{-1}(B)$, there exists an open set $C\subset
Y$ containing $B$ and such that $f^{-1}(C)\subset A$.

\item[{\rm \textbf{(ii)}}] $f$ is closed
if and only if for every point $y\in Y$ and every open set
$U\subset X$ which contains $f^{-1}(y)$, there exists in $Y$ a
neighborhood $V_y$ of the point $y$ such that $f^{-1}(V_y)\subset
U$.
\end{itemize}
\end{theorem}

We now prove a crucial lemma needed for the generalization of the
``Lokal-global-Prinzip''.

\begin{lemma}
\label{for Hausdorff} Let  $X$ be a normal first countable topological space and
$V$ a locally convex topological vector space. Let $f:X\rightarrow
V$ be a continuous map that has local convexity data and satisfies
the locally fiber connected condition. Suppose $f$ is a closed map.
Then the following hold:
\begin{itemize}
\item[{\rm \textbf{(i)}}] The projection $\pi_f: X \rightarrow X
_f$ is also a closed map.

\item[{\rm \textbf{(ii)}}] The quotient  $X_f$ is a Hausdorff
topological space.
\end{itemize}
\end{lemma}

\noindent\textbf{Proof.\ \ } \noindent\textbf{(i)} Let $[x]$ be an
arbitrary point in $X_f$ and $U\subset X$ an arbitrary open set
that contains ${\pi_f}^{-1}([x])=E_x$, the connected component of
$f^{-1}(f(x))$ that contains $x$. Denote by
$F:=f^{-1}(f(x))\setminus E_x$ which is the union of all (closed)
connected components of $f^{-1}(f(x))$ minus $E_x$. We claim that
$F$ is a closed subset of $X$.  Indeed, if $z\in \overline{F}$,
then there exist a sequence $\{z_n\}_{n \in  \Bbb N}$ in $F$ which
converges to $z$. Continuity of $f$ implies that $z\in
f^{-1}(f(x))$.
If $z\in E_x$ then any neighborhood of $z$
intersects at least one other connected component of the fiber
$f^{-1}(f(x))$ which contradicts the locally fiber connected
condition. Since $z\in f^{-1}(f(x))$ it follows that $z \in F $
and hence $F$ is closed. 
The same argument as above shows that
(LFC) implies that $E_{x}$ is also closed in $X$.

Due to the normality of $X$, there exist two open sets $U_{E_x}$ and
$W$ such that $E_x\subset U_{E_x}$, $F\subset W$, and
$U_{E_x}\cap W=\varnothing$. 
After shrinking, if necessary, we can
take $U_{E_x}\subset U$. By the closedness of $f$ there exists an
open neighborhood  $V_{f (x)}$ of  $f (x)  $ in $V$ such that
$E_x\subset f^{-1}(f(x))\subset f^{-1}(V_{f (x)})\subset
U_{E_x}\cup W$, by Theorem~\ref{caracterizare inchise} {\bf (ii)}.

The set $A:=U_{E_x}\cap f^{-1}(V_{f (x)})$ is a nonempty open
subset of $X$ which is also saturated with respect to the
equivalence relation that defines $\pi_f$ or, otherwise stated,
${\pi_f}^{-1}(\pi_f(A))=A$. Indeed, if a connected component of a
fiber of $f$ from $f^{-1}(V_{f (x)})$ intersects $U_{E_x}$, respectively $W$, then it is entirely contained 
either in
$U_{E_x}$ or in $W$ since $U_{E_x}\cap W=\varnothing$. We have
hence proved that $\pi_f(A)$ is an open neighborhood of $[x]$
since $\pi_f$ is an open map and ${\pi_f}^{-1}(\pi_f(A))\subset U$.
This shows, via Theorem \ref {caracterizare inchise} {\bf  (ii)},
that $\pi_f$ is a closed map.

\medskip

\noindent\textbf{(ii)} Since $\pi_f$ is open, it suffices to prove
that the equivalence relation that defines $X_f$ has a closed
graph in order to show that $X_f$ is Hausdorff.

For every $x\in X$ there exist a neighborhood $U_x$ that satisfies
(VN), (SLO), and (LFC). Notice that the open set $\pi _f^{-1}(\pi _f(U _x'))$ also
satisfies the (LFC) condition. By the normality of
$X$ there also exists  a neighborhood $U'_x$ of $x$ with $\overline{U'_x}\subset
U_x$.  We
shall prove that $\overline{{\pi_f}^{-1}(\pi_f(U'_x))}\subset
{\pi_f}^{-1}(\pi_f(U_x))$ which shows that for every connected
component $E_x$ of a fiber there exists a saturated neighborhood of
it which contains a smaller saturated neighborhood whose closure
still satisfies (LFC). In order to prove the above inclusion
observe that, since $\pi_f$ is closed, we have $\overline
{\pi_f(U'_x)}=\pi_f(\overline{U'_x})\subset \pi_f(U_x)$.
By the
continuity of $\pi_f$ we obtain the inclusion
$\overline{{\pi_f}^{-1}(\pi_f(U'_x))}\subset
{\pi_f}^{-1}(\overline{\pi_f(U'_x)})\subset
{\pi_f}^{-1}(\pi_f(U_x))$.

We now prove the closedness of the graph of the equivalence
relation that defines $X _f  $. Take $\{x_n\}_{n\in\mathbb{N}}$
and $\{y_n\}_{n\in \mathbb{N}}$ two convergent sequences in $X$
such that $x_n$ and $y_n$ are in the same equivalence class for
all $n\in \mathbb{N}$. Suppose that $x_n\rightarrow x$ and
$y_n\rightarrow y$. The continuity of $f$ guarantees that $f  (x)
= f  (y) $. Additionally, there exists $n_0\in \mathbb{N}$ such
that for $n>n_0$ all $x_n\in {\pi_f}^{-1}(\pi_f(U'_x))$.
Consequently $y_n\in {\pi_f}^{-1}(\pi_f(U'_x))$ since $x_n$ and
$y_n$ are in the same equivalence class and
${\pi_f}^{-1}(\pi_f(U'_x))$ is saturated. Since $y\in
\overline{{\pi_f}^{-1}(\pi_f(U'_x))}$ and
$\overline{{\pi_f}^{-1}(\pi_f(U'_x))}$ satisfies (LFC)  we obtain
that $x$ and $y$ sit necessarily on the same connected component
of the fiber $f ^{-1}  (f  (x))$ and are hence equivalent. This
shows that the graph of the equivalence relation is closed, as
required. \quad $\blacksquare$

\begin{remark}
\label{topologica} 
\normalfont Let $f:X\rightarrow V$ be a
continuous map such that $f$ has local convexity data and is
locally fiber connected. Then for every point $[x]\in X_{f}$ there
is a neighborhood $\widetilde{U}_{[x]}$ of $[x]$ such that
$\widetilde{f}:\widetilde{U}_{[x]}\rightarrow C_{x,f(x)}$ is a
homeomorphism onto its open image. Eventually, after shrinking
$\widetilde{U}_{[x]}$, we can suppose that its image is convex.
\end{remark}

Note that $X _f $ is connected since $X$  is connected.
Additionally, this remark implies that
$X_f$ is locally path connected. Therefore, $X_f $ is path
connected.
\medskip

From now on we take $V$ to be a Banach space. Define a
distance $d$ on $X_f $ in the following way: for  $[x],[y]
\in X _f $ let
$d([x],[y])$ be the infimum of all the lengths
$l(\widetilde{f}\circ \gamma )$, where $\gamma$ is a continuous curve in
$X _f $ that connects $[x]$ and $[y]$. The length is
calculated with respect to the distance $d_{V}$ defined by the
norm  on $V$. From the definition it follows that
$d_{V}(\widetilde{f}([x]),\widetilde{f}([y]))\leq d([x],[y])$ and,
by Remark \ref{topologica}, equality holds  for $[x]$ and
$[y]$ sufficiently close. Indeed, we shall prove that for any two points $[x],[y]\in
X _f $ that are both contained in an open subset $U $ that satisfies (VN), (SLO),
and (LFC), the equality $d_{V}(\widetilde{f}([x]),\widetilde{f}([y]))= d([x],[y])$
is satisfied. Due to the (VN) condition, the set $\widetilde{f}(U)$ is convex in $V$
and hence there exists a straight line $c_0$ joining $\widetilde{f}([x])$ and
$\widetilde{f}([y])$. Since $U$ satisfies (SLO) and (LFC), the map
$\widetilde{f}|_{U}:U \rightarrow  \widetilde{f}(U)$ is a homeomorphism and hence
the curve $\gamma _0:= \widetilde{f}^{-1}\circ c_0$ is continuous and joins
$[x]$ and $[y]$. Consequently,
\[
d([x],[y])= \inf_{\gamma\in \Gamma_{[x],[y]}}\{l(\widetilde{f}\circ \gamma)\}=l(c
_0)=d_{V}(\widetilde{f}([x]),\widetilde{f}([y])),
\]
where $\Gamma_{[x],[y]} :=\{\gamma:[a,b]\rightarrow X _f\mid
\text{$\gamma$ is continuous and }\gamma(a)=[x],\, \gamma(b)=[y]\}$. 

\begin{proposition}
\label{metrica} Assume the hypotheses of Lemma \ref{for Hausdorff}
and suppose that $V$ is a Banach space. The distance on $X _f $
introduced above defines a metric topology on $X _f $  that
coincides with the quotient topology of $X _f$.
\end{proposition}

\noindent\textbf{Proof.\ \ } The proof mimics those in
\cite{dazord} or \cite{hilnebplank}. The symmetry of
$d$ and the triangle inequality for $d$ are obvious. It remains to
be proved that if $d([x],[y])=0$ then $[x]=[y]$. By
contradiction, suppose that
$d([x],[y])=0$ with
$[x]\neq [y]$. Then
$d_{V}(\widetilde{f}([x]),\widetilde{f}([y]))=0
$  and hence $\widetilde{f}([x])=\widetilde{f}([y])$. Since
$X_f$ is Hausdorff we can find two disjoint open neighborhoods
$\widetilde {U}_{[x]}$ and $\widetilde {U}_{[y]}$ of  $[x] $ and
$[y] $, respectively, that behave as in Remark
\ref{topologica}. Take $r>0 $ such that the open disk
$D_r(f([x])):=\{v \in  V\mid d _V(v, f([x]))< r \} $ satisfies
 $D_r(f([x]))\cap
C_{x,f(x)}\subset \widetilde{f}(\widetilde{U}_{[x]})$. Then, for
all the curves $\gamma$ that connect $[x]$ and $[y]$ we have that
$l(\widetilde{f}\circ \gamma)>2r$. So, in order to have
$l(\widetilde{f}\circ \gamma)<2r$, $[y]$ must be contained in
$\widetilde {U}_{[x]}$, which is a contradiction since $\widetilde
{U}_{[x]}\cap\widetilde {U}_{[y]}=\varnothing$.

The construction of $d$ and
Remark \ref{topologica} show that $\widetilde{f} :X_f \rightarrow
V $ is a local isometry on the image $\widetilde{f}(\widetilde{U}_{[x]})$ of all the
small enough neighborhoods $\widetilde{U}_{[x]} $ of any point $[x]\in X _f$. It
follows that the metric topology on  $ X _f  $ coincides with the quotient topology
of $ X _f
$.\quad
$\blacksquare$

\medskip

The proof of the following lemma can be found in \cite{engelking}.

\begin{lemma}
\label{frontier} (Va\v{\i}n\v{s}te\v{\i}n)
If
$f:X\rightarrow Y$ is a closed mapping from a metrizable space $X$
onto a metrizable space
$Y$, then for every $y\in Y$ the boundary ${\rm bd
}(f^{-1}(y)):=\overline{f^{-1}(y)}\cap\overline{ \left(X \setminus
f^{-1}(y)\right)}$ is compact.
\end{lemma}

\begin{definition}
Let $X$ be a Hausdorff topological space and $f:X\rightarrow Y$ be
a continuous map. We call $f$ a \textbf{proper map} if $f$ is
closed and all fibers $f^{-1}(y)$ are compact subsets of $X$.
\end{definition}

\begin{theorem}[\cite{engelking}]
If $f:X\rightarrow Y$ is a proper map, then for every compact
subset $Z\subset Y$ the inverse image $f^{-1}(Z)$ is compact.
\end{theorem}

A converse of this theorem is available when $Y$
is a $k$-space (i.e. $Y$ is a Hausdorff topological space that is
the image of a locally compact space under a quotient mapping).

\paragraph{The finite dimensional case.} For the  next
considerations we need that $V$ is a finite dimensional vector
space with a fixed chosen inner product that defines the distance
$d_V$. The  infinite dimensional case will be discussed later.

\begin{proposition}
\label{compact bole} In the hypotheses of Lemma \ref{for
Hausdorff} with $V$ a finite dimensional Euclidean vector space we
have that $\widetilde{f}:X_f\rightarrow V$ is a proper map and
$B_r([x]):=\{[y]\in  X _f\mid d([x],[y])\leq r\}$ is compact.
\end{proposition}

\noindent\textbf{Proof.\ \ } The locally fiber connectedness
condition on $f$ and the openness of $\pi_f$ imply that
${\widetilde{f}}^{-1}(v) $ is a collection of isolated points in
$X _f $  and hence ${\widetilde{f}}^{-1}(v)={\rm bd}
\left({\widetilde{f}}^{-1}(v)\right)$, for every $v\in f(X)$. As a
consequence of the Va\v{\i}n\v{s}te\v{\i}n Lemma we obtain that
the fibers of $\widetilde{f}$ are all compact. Also,
$\widetilde{f}$ is closed since $f$ is closed and hence
$\widetilde{f}$ is a proper map.

The set $B_r([x])$ is closed in $X_f$ and, by the definition of the
distance $d$, we have that $\widetilde{f}(B_r([x]))\subset
B_r(\widetilde{f}([x]))\subset V$. Since the ball
$B_r(\widetilde{f}([x]))\subset V$ is closed and bounded it is
necessarily compact in the finite dimensional vector space $V$.
Properness of $\widetilde{f}$ implies that
${\widetilde{f}}^{-1}(B_r(\widetilde{f}([x])))$ is compact in $X
_f $. Since $B_r([x])$ is a closed subset of
${\widetilde{f}}^{-1}(B_r(\widetilde{f}([x])))$, it is necessarily
compact in $X_f$.\quad $\blacksquare$

\medskip

We can now prove one of the main results of this section which is
the generalization of the ``Lokal-global-Prinzip'' in
\cite{dazord} or \cite{hilnebplank} to the case of closed maps
that are not necessarily proper.

\begin{theorem}
\label{localglobal}
Let
$f:X\rightarrow V$ be a closed map with values in a finite
dimensional Euclidean vector space $V$ and $X$  a connected,
locally connected, first countable, and normal topological space. Assume that $f$
has local convexity data and is locally fiber connected. Then:
\begin{description}
\item[(i)] All the fibers of $f$ are connected.
\item[(ii)]$f$ is open onto its image.
\item[(iii)] The image $f(X)$ is a closed convex
set.
\end{description}
\end{theorem}

\noindent\textbf{Proof.\ \ } We begin with the argument
of~\cite{dazord, hilnebplank}. Let $[x]_{0},[x]_{1}\in X_{f}$ be
two arbitrary points and $c:=d([x]_{0},[x]_{1})$. By the
definition of $d$, we have that for every $n\in \Bbb{N}$ there
exists a curve $\gamma _{n}$ defined on the interval $[a,b]$,
connecting $[x]_{0}$ and $[x]_{1}$, and satisfying $l(\widetilde{f}\circ
\gamma _{n})\leq c+\frac{1}{n}$. Also, for every $n\in \Bbb{N}$,
let $v_{n}=(\widetilde{f}\circ \gamma _{n})(t _0)$ be the point on
the curve $\widetilde{f}\circ \gamma _{n}$ such that
$l(\widetilde{f}\circ \gamma _{n}|_{[a, t  _0]})=
\frac{1}{2}l(\widetilde{f}\circ \gamma _{n}) $.
Then there exists
a finite set of points $\{[x]_{1}^{n},...,[x]_{k_{n}}^{n}\} $  in
$X _f $ such that $\widetilde{f}^{-1}(v_{n})\cap {\rm
range}(\gamma_{n})=\{[x]_{1}^{n},...,[x]_{k_{n}}^{n}\} \subset
B_{c+1}([x]_{0})$, where $B_{c+1}([x]_{0}):=\{[x]\in X_{f}\mid
d([x]_{0},[x])\leq c+1\}$ is compact by Proposition \ref{compact
bole}. Relabeling the elements of the set $\underset{n\in
\Bbb{N}}{\bigcup }\{[x]_{1}^{n},...,[x]_{k_{n}}^{n}\}$ we obtain a
sequence included in the compact set $B_{c+1}([x]_{0})$ and,
consequently, it will have an accumulation point denoted by
$[x]_{\frac{1}{2}}$.

The definition of $d$ implies that $d([x]_{0},[x]_{\frac{1}{2}})
=d([x]_{\frac{1}{2}},[x]_{1})= \frac{c}{2}$. Repeating this
process for the pair of points $([x]_{0},[x]_{\frac{1}{2}})$ and
$([x]_{\frac{1}{2}}, [x]_{1})$ we obtain the points
$[x]_{\frac{1}{4}}$ and $[x]_{\frac{3}{4}}$ satisfying
$d([x]_{0},[x]_{\frac{1}{4}})=
d([x]_{\frac{1}{4}},[x]_{\frac{1}{2}})
=d([x]_{\frac{1}{2}},[x]_{\frac{3}{4}})
=d([x]_{\frac{3}{4}},[x]_{1})=\frac{c}{4}$. Inductively, we obtain
points $[x]_{n/2^{m}}$,  $[x]_{n'/2^{m'}}$ for $0\leq n\leq
2^{m}$,  $0\leq n'\leq 2^{m'}$, such that
\begin{equation}
\label{modu} d([x]_{n/2^{m}},[x]_{n^{^{\prime }}/2^{m^{^{\prime
}}}}) =c\left| \frac{n}{2^{m}}-\frac{n^{^{\prime }}}{2^{m^{\prime
}}}\right|.
\end{equation}

We can extend the map $n/2^{m}\mapsto [x]_{n/2^{m}}$ to a
continuous map $\gamma :[0,1]\rightarrow X_{f}$ such that
\begin{equation}
\label{modul}
d(\gamma(t),\gamma (t^{^{\prime }})) =c|
t-t^{^{\prime }}|.
\end{equation}
To see this, note that every $t\in [0,1]$ can be approximated by a
sequence of the type $n_{k}/2^{m_{k}}$. The corresponding points
$[x]_{n_{k}/2^{m_{k}}}$ are contained in the compact set
$B_{c+1}([x]_{0})$ and hence they have an accumulation point
$[x]_{t}$.
It is now easy to see, using~(\ref{modu}), that
$[x]_{t}$ does not depend on the sequence $n_{k}/2^{m_{k}}$ and
that the curve $\gamma$ constructed in this way is continuous.

Remark~\ref{topologica} and~(\ref{modul}) imply  that,
locally,  $d_{V}((\widetilde{f}\circ \gamma)
(t),(\widetilde{f}\circ \gamma)(t^{^{\prime }})) =c\left|
t-t^{^{\prime }}\right| $ which shows that
$\widetilde{f}\circ \gamma $ is locally a straight line. Due
to~(\ref{modul}),
$\widetilde{f}\circ \gamma $ is necessarily a straight line that
goes through $\widetilde{f}([x]_{0})$ and
$\widetilde{f}([x]_{1})$.
This proves the convexity of $f (X)
$. Since $f $ is a closed map the set $f(X)$ is closed in $V $
which proves {\bf (iii)}.

In order to prove the connectedness of the fibers of $f$ let
$[x]_{0},[x]_{1}\in X_{f}$ be two arbitrary points such that
$v:=\widetilde{f}([x]_{0})=\widetilde {f}([x]_{1})$ and
$c:=d([x]_{0},[x]_{1})$. Any curve that connects these two points
is mapped by $\widetilde{f}$ into a loop based at $v$. We shall
prove that $c=0$, which implies that $[x]_{0}=[x]_{1}$ and hence
that the fibers of  $f$ are connected.
Let $\gamma $ be the curve constructed above. Then the range of
$\widetilde{f}\circ \gamma $ is a segment that contains $v$.
We will prove by contradiction that this segment consists of just
one point which is $v$ itself. 

Suppose that this is not true. Since
$\widetilde{f}\circ \gamma $ is a loop based at $v$ and at the same time a straight
line, there exists a turning point $v_{0}:=(\widetilde{f}\circ \gamma)
(t_{0})$ on the segment $\widetilde{f}\circ \gamma $ such that for
$t\leq t_{0}$ we approach $v_{0}$ and for $t^{^{\prime}}\geq
t_{0}$ we  move away from $v_{0}$ staying on the same segment
which is the range of $\widetilde{f} \circ \gamma$. Otherwise
stated, ${\rm range}(\widetilde{f}\circ\gamma|_{[t,t_{0}]})= {\rm
range}(\widetilde{f} \circ\gamma|_{[t_{0},t^{^{\prime}}]})$ and
hence in a neighborhood of $\gamma (t _0)$ the map $\widetilde{f}
$ is not injective. However, since $f$  is locally fiber connected
the map $\widetilde{f} $ is locally injective, which is a
contradiction. This proves {\bf (i)}. Note that from
$c:=d([x]_{0},[x]_{1})$ and
$d_{V}(\widetilde{f}([x]_{0}),\widetilde{f}([x]_{1}))= d _V(v,v)
=0$ we cannot conclude that $c=0 $ since the equality between the
two metrics holds only locally.

The openness of $f$ is implied by Corollary \ref{imediat}.\quad
$\blacksquare$

\begin{remark}
\label{restriction}
\normalfont 
This result remains true if we replace the vector space $V$ by
a convex subset $C$ of $V$. More specifically, Theorem \ref{localglobal} is
valid when we apply it to a map $f:X \rightarrow C$, with $C$ a convex subset
of $V$. In particular, this allows us to generalize to the case of closed maps
 Theorem 4.2 of Lerman et al. \cite{LMTW} and  Theorem
3.3 of Weinstein
\cite{Weinstein}, initially stated for proper maps. This remark is also important
later on when we generalize several other classical convexity theorems.
\end{remark}

\paragraph{The infinite dimensional case.} Now let $V$ be an
infinite dimensional Banach space and $d_V$ the distance induced
by the norm. Analyzing  the  proof of Proposition \ref {compact
bole} it can be seen that the compactness of the closed balls in
the norm topology was essential. Thus, there is no direct analog
of this statement in infinite dimensions. What is needed is a
second topology on $V$ whose closed balls are compact. A natural
hypothesis is  that $V $ is the topological dual space of another
Banach space $W $ for then Alaoglu's Theorem guarantees that $B
_r(0):=\{v \in V \mid \|v\|\leq r\} $ is weak$^*$ compact. We
shall denote by $(V, \|\cdot\|) $ the Banach space $V $ endowed
with the norm topology (and hence this is the Banach space dual
$W^* $) and by $(V, w^*)$ the space $V $ endowed with the weak$^*$
topology of $W ^*$. Since the weak$^*$ topology is weaker than
then norm topology we have
\begin{itemize}
\item if $f:X \rightarrow (V, \|\cdot \|) $ is continuous
then $f:X \rightarrow (V, w^*) $ is continuous;
\item if $f:X \rightarrow (V, w^*) $ is closed then $f:X
\rightarrow (V, \|\cdot \|) $ is closed.
\end{itemize}

The analog of Proposition \ref {compact bole} is the following.

\begin{proposition}
\label{compact inf ball}
Assume the hypotheses of Lemma \ref{for Hausdorff} with $V = W^*$,
for $W$ a Banach space. If $f:X \rightarrow (V, \|\cdot \|) $ is
continuous and $f:X \rightarrow (V, w^*) $ is closed then
$\widetilde{f}:X_f\rightarrow V$ is a proper map relative to both
topologies on $V$ and $B_r([x]):=\{[y]\in  X _f\mid d([x],[y])\leq
r\}$ is compact in $X_f $. 
\end{proposition}

\noindent\textbf{Proof.\ \ } The locally fiber connectedness
condition for $f$ and the openness of
$\pi_f$ imply that ${\widetilde{f}}^{-1}(v) $ is a collection of
isolated points in $X _f $  and hence
${\widetilde{f}}^{-1}(v)={\rm bd}
\left({\widetilde{f}}^{-1}(v)\right)$, for every
$v\in f(X)$.
Since $\widetilde{f}:X_f \rightarrow (V, w^*) $ is a closed map
then so is
$\widetilde{f}:X_f\rightarrow (V, \| \cdot \|)$.
By Va\v{\i}n\v{s}te\v{\i}n's Lemma
the fibers of $\widetilde{f}$ are all compact. Also
$\widetilde{f}:X_f \rightarrow (V, w^*) $ is closed and hence
$\widetilde{f}:X_f \rightarrow (V, w^*) $ is a proper map.

The set $B_r([x])$ is closed in $X_f$ and by the definition of the
distance $d$ we have that $\widetilde{f}(B_r([x]))\subset
B_r(\widetilde{f}([x]))\subset V$. The ball
$B_r(\widetilde{f}([x]))\subset V$ is weak$ ^\ast $ compact and
since $\widetilde{f}:X_f \rightarrow (V, w^*) $ is a proper map we
have that ${\widetilde{f}}^{-1}(B_r(\widetilde{f}([x])))$ is
compact in $X _f $. As $B_r([x])$ is a closed subset of
${\widetilde{f}}^{-1}(B_r(\widetilde{f}([x])))$, it is necessarily
compact in $X_f$.\quad $\blacksquare$

\medskip

Using Proposition~\ref{compact inf ball}, one can
generalize Theorem~\ref{localglobal} to the context of
maps with an infinite dimensional range by virtually copying its proof.

\begin{theorem}
\label{localglobal infinite}
Let $(V,\| \cdot \|)$ be a Banach space such that $V = W^*$,
for $W$ a Banach space. Let $f:X \rightarrow (V, \|\cdot \|) $ be
a continuous map  and $f:X \rightarrow (V, w^*) $
closed, where $X $ is a connected,
locally connected, first countable, and normal topological space. Assume that $f$
has local convexity data and is locally fiber connected. Then:
\begin{description}
\item[(i)] All the fibers of $f$ are connected.
\item[(ii)]$f:X \rightarrow  (V, w ^\ast )$ is open onto its image.
\item[(iii)] The image $f(X) \subset (V,w ^\ast)$ is a closed
convex set.
\end{description}
\end{theorem}

\begin{remark}
\normalfont Since the weak$ ^\ast $ topology is weaker than the
norm topology the previous theorem also implies that $f:X
\rightarrow  (V, \| \cdot \|)$ is open onto its image and that $f (X)  $ is
closed in $(V, \| \cdot \|)$ (and obviously also convex since the notion of
convexity is not related to the topology). Moreover, as the fibers of 
$f$ are connected the map $ \widetilde{f}:X _f \rightarrow f (X) $ is a
homeomorphism when $f (X) $ is considered as a topological
subspace of either $(V, \| \cdot \|) $ or  $(V, w ^\ast )$. This
apparently strong conclusion is actually  guaranteed to hold
automatically for an important class of normed spaces $(V, \|\cdot
\|)$. For instance, as $V= W ^\ast $, if $W$ is reflexive then the
weak and the weak$ ^\ast  $ topologies in $V$ coincide. Moreover,
Mazur's Theorem~\cite{mazur} guarantees that the weak  and the
norm closures of a convex set on a normed space coincide. The same  interplay
between the weak$^\ast $ and norm topology can be found in the infinite
dimensional convexity  results of Bloch-Flaschka-Ratiu~\cite{bfr} and
Neumann~\cite{neumann}.
\end{remark}

\section{Openness and local convexity for momentum maps}
\label{Openness and local convexity for momentum maps}

In the previous section we presented sufficient conditions for the
convexity of the image of a map that has both local convexity data and is locally fiber connected. More specifically, in
Theorem~\ref{conv} we saw that if the map is open onto its image,
then the image is locally convex and therefore convex if it is
also closed. In this section we will characterize the situations
in which the momentum map associated to a compact and symplectic
Lie group action on a symplectic manifold is an open map onto its
image. We will also give several generalizations of convexity
results found in the literature. Throughout this paper, all
manifolds are assumed to be paracompact, that is, they are Hausdorff spaces and
every open cover has a locally finite open refinement.

\begin{definition}
Let $V$ be finite dimensional vector space.
\begin{description}
\item[{\rm \textbf{(i)}}] A subset $K\subset V$ is called {\bfi
polyhedral} if it is the intersection of a finite family of closed
halfspaces of $V$. Consequently, a polyhedral subset of $V$ is
closed and convex.
\item[{\rm \textbf{(ii)}}] A subset $K\subset
V$ is called {\bfi locally polyhedral} if for every $x\in K$ there
exists a polytope $P_{x}$ in $V$ such that $x\in \operatorname{int}(P_{x})$ and
$K\cap P_{x}$ is a polytope.
\end{description}
\end{definition}

In order to proceed we need further preparation. First we
specialize the Marle-Guillemin-Sternberg Normal Form Theorem to
the case of  torus actions in the formulation of
\cite{hilnebplank}; see also
\cite{marle, GuSt1984, OrRa2004} for general normal form theorems
and  their proofs.

\begin{theorem}
\label{MGS theorem}
Let $(M,\omega )$ be a symplectic manifold and let $T$ be a
torus acting on $M$ in a globally Hamiltonian
fashion with invariant momentum map $J_{T}:M\rightarrow
\mathfrak{t}^{\ast }.$ Let $m\in M$ and $T_{0}=(T_{m})^{0}$ be
the connected component of the stabilizer $T_{m}$. Let $T _1\subset
T$  be a subtorus such that
$T=T_{0}\times T_{1}$.
Then:
\begin{itemize}
\item[{\rm \textbf{(i)}}] There exist a symplectic vector space
$(V,\omega _{V})$, a $T$-invariant open neighborhood  $U\subset M$
of the orbit $T \cdot  m $, and a symplectic  covering of   a
$T$-invariant open subset   $U'$ of  $T _1 \times
\mathfrak{t}_1^\ast \times  V $ onto $U$  under which the $T
$-action on $M$  is modeled by
\begin{eqnarray*}
(T_{0}\times T_{1})\times ((T_{1}\times
\mathfrak{t}_{1}^{\ast })\times V)
&\rightarrow &((T_{1}\times \mathfrak{t}_{1}^{\ast })\times V) \\
((t_{0},t_{1}),(t_{1}^{^{\prime }},\beta ,v)) &\mapsto
&(t_{1}t_{1}^{^{\prime }},\beta ,\pi (t_{0})v),
\end{eqnarray*}
where $\pi :T_{0}\rightarrow Sp(V)$ is a symplectic
representation.

\item[{\rm \textbf{(ii)}}] There exists a complex structure $I$
on $V$ such that $\langle v,w\rangle:=\omega_{V}(Iv,w)$
defines a positive definite scalar product on $V$. Then
$V=\oplus _{\alpha \in \mathcal{P}_{V}}V_{\alpha }$, where
$V_{\alpha}:=\{v\in V\mid Y\cdot v=\alpha (Y)Iv,\text{for
all} \; Y\in \mathfrak{t}_{0} \}$ and
$\mathcal{P}_{V}:=\{\alpha \in
\mathfrak{t}_{0}^{\ast } \mid V_{\alpha }\neq \{0\}\}$. The
corresponding $T$-momentum map $\Phi:T^{\ast }(T_{1})\times
V\rightarrow \mathfrak{t}_{1}^{\ast}\times
\mathfrak{t}_{0}^{\ast }\simeq \mathfrak{t}^{\ast}$ is given by
\[
\Phi\left((t_{1},\beta ),\sum_\alpha v_\alpha \right)  =  \Phi
(1,0.0)+\left(\beta ,\frac{1}{2}\sum\limits_{_{\alpha \in
\mathcal{P}_{V}}}|| v_{\alpha }||^{2} \alpha\right).
\]
\end{itemize}
\end{theorem}

\medskip

Notice that the original version of the ÊMarle-Guillemin-Sternberg Normal Form
Theorem provides the twisted product $(T _0 \times  T _1)\times _{T
_0}(\mathfrak{t}_1 ^\ast  \times  V )$ as a $T$-invariant local model for $M$. This
is equivariantly diffeomorphic to $T _1\times  \mathfrak{t}_1^\ast \times  V $
 via the map
\begin{equation*}
\begin{array}{ccc}
(T _0 \times  T _1)\times _{T_0}(\mathfrak{t}_1 ^\ast  \times  V )& \longrightarrow
& T _1\times  \mathfrak{t}_1^\ast \times  V \\
 \left[(t_{1},t_{0}), \eta, v \right]&\longmapsto & (t _1, \eta, t _0 \cdot  v).
\end{array}
\end{equation*}

The next result shows that the momentum maps of globally
Hamiltonian torus actions always have local convexity data with
closed cones and are locally fiber connected. In fact, the
associated cones are closed in $\mathfrak{t}^{\ast }$.

\begin{theorem}
\label{conuri} Let $(M,\omega )$ be a symplectic manifold and let
$T$ be a torus acting on $M$ in a globally Hamiltonian fashion
with invariant momentum map $J_{T}:M\rightarrow \mathfrak{t}^{\ast
}$. Then there exist an arbitrarily small neighborhood $U$ of $m$
and a convex polyhedral cone $C_{J(m)}\subset \mathfrak{t}^{\ast
}$ with vertex $J_{T}(m)$ such that:
\begin{description}
\item[(i)]$J_{T}(U)\subset C_{J_{T}(m)}$ is an
open neighborhood of $J_{T}(m)$ in $C_{J_{T}(m)}$;

\item[(ii)] $J_{T}:U\rightarrow C_{J_{T}(m)}$ is
an open map;

\item[(iii)] If $\mathfrak{t}_{0}$ is the Lie
algebra of the stabilizer $T_{m}$ of $m$, then
$C_{J_{T}(m)}=J_{T}(m)+\mathfrak{t}_{0}^{\perp}+
\operatorname{cone}(\mathcal{P}_{V})$;

\item[(iv)] $J_T ^{-1}(J_T(m)) \cap U$ is
connected for  all $m \in U $.
\end{description}
\end{theorem}

Here is a sketch of the proof; for details see
\cite{hilnebplank}. Recall that
$\operatorname{cone}(\mathcal{P}_V) : = \{ \sum_j a_j \alpha_j
\mid a_j \geq 0 \}$.  By Theorem
\ref{MGS theorem}, it suffices to work with the momentum map
$\Phi$. For small neighborhoods
$B_{\mathfrak{t}_{1}^{\ast }}$ and $B_{V}$ of the origin in
$\mathfrak{t}_{1}^{\ast }$ and $V$ respectively, the restriction
of $\Phi$ to $U : = T_{1}\times B_{\mathfrak{t}_{1}^{\ast }}
\times B_{V}$ takes values in the polyhedral closed convex cone
$C_{J_{T}(m)}=J_{T}(m)+ \mathfrak{t}_{1}^{\ast} +
\operatorname{cone}(\mathcal {P}_{V})$, where
$\operatorname{cone}(\mathcal {P}_{V})$ denotes the cone
generated by the finite set
$\mathcal{P}_{V}=\{\alpha_{1},...,\alpha _{n}\}$ of
$T_0$-weights for the action on $V$;
$\operatorname{cone}(\mathcal {P}_{V})$ is clearly closed. In order
to prove that $J_{T}$ satisfies the conditions of the theorem, we
decompose it in two maps
$\varphi_{1}: (t_{1},\beta ,\sum_\alpha v_\alpha)\mapsto (\beta
,(\|v_{\alpha }\| ^{2}))$ and
$\varphi _{2}:(\beta,a_{1},...,a_{n})\mapsto
(\beta,\frac{1}{2}\sum_j a_{j}\alpha _{j})$. We have
$J_{T}=\varphi _{2}\circ \varphi _{1}+J_{T}(m)$. One proves that
$\varphi_1, \varphi_2 $ are open onto their images
and have connected fibers, so $J_T $ has the same properties.
\medskip

We now state a generalization of the Atiyah-Guillemin-Sternberg Convexity
Theorem for non-compact manifolds.

\begin{corollary}
Let $M$ be a paracompact connected  symplectic manifold on which a
torus $T$ acts in a Hamiltonian fashion. Let  $J_T:M\rightarrow \mathfrak{t}^{*}$ be
an associated momentum map which we suppose is closed. Then the image $J_T(M)$
is a closed convex locally polyhedral subset in
$\mathfrak{t}^\ast$. The fibers of $J_T$ are connected and  $J_T $
is open onto its image.
\end{corollary}

\noindent\textbf{Proof.\ \ } This is a consequence of Theorems
\ref{conuri} and \ref{localglobal}. The image is locally
polyhedral since $J_T $ is open onto its image and the associated
cones are polyhedral. \quad $\blacksquare$
\medskip

Our approach to convexity also allows us to generalize a result
due to Prato~\cite{prato}.

\begin{theorem}
\label{prato theorem better}
Let $M$ be a connected symplectic manifold on which a
torus $T$ acts in a Hamiltonian fashion with associated invariant
momentum map $J_T:M\rightarrow \mathfrak{t}^{*}$.
\begin{description}
\item [(i)] If there exists $\xi\in \mathfrak{t} $ such that the map
$J_T ^\xi\in  C^\infty(M)$ defined by $J _T ^\xi:=\langle J_T, \xi \rangle$  is
proper, then the image
$J_T(M)$ is a closed convex locally polyhedral subset in
$\mathfrak{t}^\ast$. Moreover, the fibers of $J_T$ are connected and  $J_T $
is open onto its image.
\item  [(ii)] If there exists an integral element $\xi \in  \mathfrak{t}$
such that $J_T^\xi  $ is a proper function having a minimum as its
unique critical value then $J _T (M) $ is the convex hull of a
finite number of affine rays in $\mathfrak{t}^\ast $ stemming from the
images of $T$-fixed points.
\end{description}
\end{theorem}

\noindent\textbf{Proof.\ \ } {\bf (i)} We will carry out the proof using our
generalization of the ``Lokal-global-Prinzip'' (Theorem~\ref{localglobal}).
According to this result all that needs to be proved is that in
the presence of our hypotheses $J_T$ is a closed map. This will be
shown by verifying that
$\overline{ J _T (A)} \subset J _T (\overline{A})$, for
any subset $A\subset M $.

We start by noticing that the map $J_T^\xi \in
C^\infty(M) $ can be written as $J_T^\xi=b \circ \pi
\circ J_T$, where $\pi: \mathfrak{t}^\ast \rightarrow
{\rm span}\{ \xi\} ^\ast $ is the dual of the inclusion ${\rm
span}\{ \xi\} \hookrightarrow \mathfrak{t}$ and $b: {\rm span}\{
\xi\} ^\ast  \rightarrow \mathbb{R} $ is the linear isomorphism
obtained as the map that assigns to each element in ${\rm span}\{
\xi\} ^\ast $ its coordinate in the dual basis of $\{\xi \}$ as a
basis of ${\rm span}\{ \xi\} $. Let $ \mu \in \overline{
J _T (A)} $ be arbitrary and $\{\mu _n \}_{n \in \mathbb{N}}
\subset J_T (A) $ a sequence such that  $\mu
_n\rightarrow  \mu $. Let $\{x _n \}_{n \in  \mathbb{N}} \subset A
$  be a sequence such that  $ J_{T}(x _n)= \mu _n $. By
continuity we have that $J ^\xi_T (x _n)=b \circ \pi
\circ J_T(x _n) \rightarrow b \circ  \pi (\mu) $. Since
by hypothesis $J _T ^\xi $ is a proper map there exists a
convergent subsequence $x_{n _k} \rightarrow x \in  \overline{A} $
and hence $J_T(x_{n _k} ) \rightarrow J _T (x) =
\mu  $, which shows that $\mu \in J_T (\overline{A})$, as
required. Part {\bf (ii)} is the original result of
Prato~\cite{prato}. \quad $\blacksquare$

\medskip

We now recall some standard notions from the theory of proper Lie
group actions. Let $M $ be a manifold and $G $  a Lie group
acting properly on it. The orbit $G\cdot m$ is called {\bfi
regular} if the dimension of nearby orbits coincides with
the dimension of $G\cdot m$. Let $M^{reg}$ denote the union of
all regular orbits. For every connected component
$M^{0}$ of $M$ the subset $M^{reg}\cap M^{0}$ is connected,
open, and dense in $M^{0}$. Note that if $U$ is a $T$-invariant
connected submanifold, then the set of regular points for the
$T$ induced action on $U$ equals $U \cap M^{reg}$. This shows
that $U \cap M^{reg}$ is open, dense, and connected in $U$.
\medskip

Next, we give necessary topological conditions on the image of
$J_{T}$ that ensure that $J_T$ is open onto its image. Later
we will show that these conditions are also sufficient. We recall
that a {\bfi  region} in a topological space is, by
definition, a connected open set.

\begin{proposition}
\label{conectregion}
Let $J_{T}:M\rightarrow \frak{t}^{*}$ be the momentum map of
a torus action on the connected symplectic manifold $(M,\omega
).$ Suppose that $J_{T}$ is open onto its image. Then the
complement
$\mathbf{C}J_{T}(M^{reg}): = J_T(M)
\setminus J_{T}(M^{reg})$ does not disconnect any region in
$J_{T}(M)$.
\end{proposition}

\noindent\textbf{Proof.\ \ }
Suppose there exists a region $V \subset J_{T}(M)$ (relative to
the induced topology from $\mathfrak{t}^\ast$) such that
$V\backslash \mathbf{C}J_{T}(M^{reg})$ is disconnected and hence
$V\backslash  \mathbf{C}J_{T}(M^{reg})= A\cup B$, where $A$ and
$B$ are open in $V $ and $A\cap B=\varnothing$. Moreover, since
$J_T(M^{reg})$ is dense in $J_T(M)$, we have that
$J_T(M^{reg})\cap V = V\backslash \mathbf{C}J_{T}(M^{reg}) =
A\cup B $ is dense in $J_T(M) \cap V = V $. Hence $\overline{(A
\cup B)} \cap V = V = (\overline{A} \cap V) \cup(\overline{B}
\cap V)$.

We claim that there exists an element $v\in
\mathbf{C}J_{T}(M^{reg})\cap V$ such that any neighborhood $V_{v}\subset V$ of $v $ in $V $ is  disconnected
by $ \mathbf{C}J_{T}(M^{reg})\cap V$.
Suppose that this claim is
false. Then for every $v \in \mathbf{C}J_{T}(M^{reg})\cap V$ there
would exist a neighborhood $V_{v}\subset V$ of $v $ such
that $V_{v}\backslash \mathbf{C}J_{T}(M^{reg})$ is connected.
Therefore, we must have that either $V_{v}\backslash
\mathbf{C}J_{T}(M^{reg})\subset A$ or $V_{v}\backslash
\mathbf{C}J_{T}(M^{reg})\subset B$. Thus, either $v \in
\overline{A}\cap V $ or $v \in \overline{B}\cap V$. Since $v \in
\mathbf{C}J_{T}(M^{reg})\cap V$ is arbitrary, this shows that
$(\overline{A} \cap V) \cap (\overline{B} \cap V) = \varnothing $.
This contradicts the
connectivity of $V $ and hence there exist a
$v\in\mathbf{C}J_{T}(M^{reg})\cap V$ such that any neighborhood $V_{v}\subset V$ is disconnected by
$\mathbf{C}J_{T}(M^{reg})$.

Take an arbitrary element $x\in J_{T}^{-1}(v)$ and $U_{x}$ a small
neighborhood of $x$ such that $J_{T}(U_{x})\subset V$ is an open
neighborhood of $v$ in $J _T(M)$; this holds because $J_T $ is
open onto its image by hypothesis. Then, by assumption,
$J_{T}(U_{x})\cap J_{T}(M^{reg})$ is disconnected. Taking the $T
$-saturation of $U_x $ we get a $T $-invariant neighborhood whose
image is in $V $ since $J_T $ is $T $-invariant. Thus, we can
assume that $U_x $ is $T$-invariant and then the set of regular
points for the induced $T$-action on $U _x$ equals the set
$U_x\cap M^{reg}$ which in turn is open, dense, and connected in
$U_x$.

Let $E: = \{z \in U_x \mid J_T(z) \in
\mathbf{C}J_T(M^{reg})\}$.  Since we can write $U_x=E\cup D$ with
$D:=U _x\setminus E $, by the construction of $E$ we have
$J _T(E)=J _T(U _x)\cap \mathbf{C}J_T(M^{reg})$ and $J _T(D)=J_{T}(U_{x})\cap
J_{T}(M^{reg})$. Now, since $E \subset U_x \setminus M^{reg}$, the
inclusion $U_x\cap M^{reg}\subset D$ also holds. Because $U_x\cap
M^{reg}$ is dense and connected in $U_x$ so is  $D$ in $U_x$.
But this is a contradiction with the fact that
$J _T(D)=J_{T}(U_{x})\cap J_{T}(M^{reg})$ is disconnected. This
proves the result. \quad $\blacksquare$

\medskip

Now we will prove the converse of Proposition \ref{conectregion}
in the case in which the momentum map has connected fibers. For this
we need a preparatory lemma.

\begin{lemma}
\label{affine space}
Let $(M, \omega) $ be a connected symplectic manifold and $J_T :
M \rightarrow \mathfrak{t}^\ast$ be the invariant momentum map
associated to a canonical $T$-action on $M$. Then
$J_T|_{M^{reg}}: M^{reg} \rightarrow J_T(M)$ is an open map.
In particular $J_{T}(M^{reg})$ is an open dense subset of
$J_{T}(M)$.
\end{lemma}

\noindent\textbf{Proof.\ \ } We shall prove that for each point in
$M^{reg}$ there is an open neighborhood such that the
restriction of $J_T $ to this neighborhood is an open map onto its
image. Let $x_0 \in M^{reg} $ be an arbitrary point. By the openness
and the $T$-invariance of $M^{reg}$ we can find an open connected $T
$-invariant neighborhood $U_{x_0}$ of $x_0 $ included in
$M^{reg}$. Therefore, for any $x \in U_{x_0}$, we have $\dim T\cdot
x =
\dim T \cdot x_0 = \dim T/T_0 = \dim T_1 $, where $T_0: = (T_{x_0})
^0$ and
$T = T_0 \times T_1$. Eventually shrinking $U_{x_0}$, using
Theorem \ref{MGS theorem}, we can work with the normal form.
Recall that the original action is symplectically and
$T$-equivariantly transformed to the action
\begin{eqnarray*}
(T_{0}\times T_{1})\times ((T_{1}\times
\mathfrak{t}_{1}^{\ast })\times V)
&\rightarrow &((T_{1}\times \mathfrak{t}_{1}^{\ast })\times V) \\
((t_{0},t_{1}),(t_{1}^{^{\prime }},\beta ,v)) &\mapsto
&(t_{1}t_{1}^{^{\prime }},\beta ,\pi (t_{0})v),
\end{eqnarray*}
where $\pi :T_{0}\rightarrow Sp(V)$ is a linear symplectic
representation. Since the isotropy subgroup of this action at
the point $(t_{1}^{^{\prime }},\beta ,v)$ equals $\{ t_0 \in T_0
\mid \pi(t_0)v = v \} \times \{e\} \subset T _0 \times  T _1$, the
condition that it be equal to
$T_0 \times \{e\} $ implies that the representation $\pi$ is
trivial. Therefore all its weights are zero. By Theorem \ref{conuri}
we conclude that
$C_{J_{T}(x_0)}=J_{T}(x_0)+ \mathfrak{t}_{1}^{\ast} +
\operatorname{cone}(\mathcal {P}_{V}) = J_{T}(x_0)+
\mathfrak{t}_{1}^{\ast}$ and that $J_T: U_{x_0} \rightarrow
C_{J_{T}(x_0)} = J_{T}(x_0)+ \mathfrak{t}_{1}^{\ast}$ is an open
map.

Note that $J_T(M^{reg}) \subset J_T(x_{0}) +
\mathfrak{t}_{1}^{\ast}$ for some (and hence any) $x_{0} \in
M^{reg} $ and  $T_1 $ is the torus whose Lie algebra is
$\mathfrak{t}_0^\perp $, where $\mathfrak{t}_0 $ is the isotropy
algebra of a regular point in $M $ and the perpendicular is taken
relative  to an a priori chosen $T $-invariant inner product on
$\mathfrak{t}$. Indeed, using the well-chained property of
$M^{reg}$ any two points in $M^{reg}$ can be linked by a finite
chain formed by the open neighborhoods constructed
above. The image of each such neighborhood lies in a translate of
$\mathfrak{t}^\ast_1 $ and since the neighborhoods intersect
pairwise, all theses affine spaces coincide. Thus, $J_T(M^{reg})$
lies in just one translate of $\mathfrak{t}_1^\ast$.  By the
density of $M^{reg} $ in $M $ and the closedness of the affine
space in $\mathfrak{t}^\ast$ it follows that $J_T(M)$ lies in the
same affine space.

Hence, we have shown that for any $x_0 \in M^{reg}$ there exists an
open neighborhood $U_{x_0} \subset M^{reg}$ such that
$J_T(U_{x_0})$ is open in a given translate of $\mathfrak{t}_1
^\ast$. Therefore, $J_T(U_{x_0}) $ is open in $J_T(M) $.
\quad $\blacksquare$

\begin{proposition}
\label{connected fibers}
Let $J_{T}:M\rightarrow \mathfrak{t}^{\ast }$ be the momentum
map of a torus action. Assume that $J_T $ has connected fibers
and that $\mathbf{C}J_{T}(M^{reg}): = J_T(M)
\setminus J_{T}(M^{reg})$ does not disconnect any
region in $J_{T}(M)$. Then $J_{T}$ is open onto its image.
\end{proposition}

\noindent\textbf{Proof.\ \ } If $M$ has more than one connected
component then the $J_{T}$-images of any two components do not
intersect, for otherwise this would contradict the connectedness
of the fibers. Since connected components of $M$ are necessarily
$T$-invariant,  we can suppose without loss of generality that $M$
is connected. We will establish the openness of $J _T  $ onto its
image through Proposition~\ref{criteriu}. In order to apply this
result it only remains to be shown that for any $v \in J _T (M) $
and any neighborhood $V _v $, the pre-image $J _T ^{-1}(V _v  )$ is
connected in $M$.

By Lemma \ref{affine space} we know that $J_{T}|_{M^{reg}}:
M^{reg} \rightarrow J_T(M)$ is an open map. Denote, as in the
Lemma \ref{affine space}, by $\mathfrak{t}_1 ^\ast$ the dual of
the subtorus whose translate contains $J_T(M)$.

Since $M $ is path connected, $J_T(M) $ is also path connected
and thus it is also locally connected. Let $J_{T}(x)\in J_{T}(M)$
be arbitrary. Choose a small neighborhood $V_{J_{T}(x)}$ of
$J_{T}(x)$ in
$\mathfrak{t}_1^{\ast }$ such that $V:=V_{J_{T}(x)}\cap J_{T}(M)$
is a region in $J_T(M) $. Then $V_{0}:=V_{J_{T}(x)}\cap
J_{T}(M^{reg})$ is connected  due to the hypothesis that
the region $V:=V_{J_{T}(x)}\cap J_{T}(M)$ cannot be
disconnected by removing $\mathbf{C}J_{T}(M^{reg})$. Now we are
in the hypotheses of Lemma  \ref{conex}. Indeed,
$J_T|_{M^{reg}}: M^{reg} \rightarrow J_T(M) $ is an open map and
we just showed that $V_0 \subset J_T(M^{reg})$ is connected. Any
fiber of $J_T $ is connected by hypothesis. Since such a fiber
is $T$-invariant, the set of its regular points for the
$T$-induced action is open dense and connected in it. If $v \in
J_T(M^{reg})$, then $J_T ^{-1}(v) \cap M^{reg}$ is connected.
Now applying Lemma \ref{conex} we
conclude that $J_{T}^{-1}(V_{0}) \cap M^{reg}$ is
connected. Since $J_{T}^{-1}(V_{0}) \cap M^{reg}$ is dense in
$J_{T}^{-1}(V_{0})$ it follows that $U_{0}:= J_{T}^{-1}(V_{0})$
is connected.

Next we will show that $U_{0}$ is dense in $U:=J_{T}^{-1}(V).$
Indeed, if this is not true, then there exist an element $x_{0}\in
U\backslash U_{0}$ and a neighborhood $U_{x_{0}}$ that does not
intersect $U_{0}.$ For the open set $U^{\prime }:=U\cap
U_{x_{0}}\cap M^{reg}\neq\varnothing$ we have that
$J_{T}(U^{\prime }) \subset V_0$ is open in $V_{0}.$ So there
exists an element $v_{0}\in J_{T}(U^{\prime })$ such that
$J_{T}^{-1}(v_{0})\cap U_{x_{0}}\neq \varnothing $ and
$J_{T}^{-1}(v_{0})\cap U_{0}\neq \varnothing$. But
$J_{T}^{-1}(v_{0})\subset U_{0}$ which contradicts the
assumption that $U_{x_{0}}\cap U_{0}=\varnothing $.

By the connectedness of $U_{0}$ and the fact that it is dense in
$U$ we obtain that $U$ is connected and hence the result follows
by Proposition \ref{criteriu}.
\quad $\blacksquare$

\medskip

We summarize Propositions \ref{conectregion} and \ref{connected
fibers} in the following theorem. Since we assume that
$J_T $ has connected fibers, we work on a possibly disconnected
symplectic manifold and we apply these propositions to each
connected component separately.

\begin{theorem}
\label{prima}
Let $J_{T}:M\rightarrow \mathfrak{t}^{\ast }$ be the momentum
map of a torus action which has connected fibers. Then
$J_{T}$ is open onto its image if and only if
$\mathbf{C}J_{T}(M^{reg})$ does not disconnect any region in
$J_{T}(M)$. Moreover, the image of the momentum
map is locally convex and locally polyhedral.
\end{theorem}

\begin{remark}
\label{topological convexity remark} \normalfont The proof of
Theorem \ref{prima} does not use in an essential way the finite
dimensionality of the manifold and the torus. A careful look at
the proof shows that same result holds for continuous maps $f: X
\rightarrow V $ that have local convexity data with closed cones,
where $X$ is a connected Hausdorff topological space and $V$ is a
locally convex topological vector space.  Additionally, assume
that
\begin{itemize}
\item $X$ is path connected;
\item there exists an open dense connected subset $X' \subset X$
such that $f|_{X'}$ is open in $f(X) $;
\item for any $v \in V $, the fiber $f ^{-1}(v) $ is connected
and $f ^{-1}(v) \cap X' $ is connected.
\end{itemize}
Then $f $ is open onto its image if and only if
$\mathbf{C}f(X'): =f(X)
\setminus f(X')$ does not disconnect any region in $f(X) $.
Moreover, the image of $f $ is locally convex. If $f(X)$ is, in
addition, closed then it is convex.
\end{remark}

We illustrate the above results with two examples, one
in which the momentum map is  open onto its image and another one
in which  it is not. This information is obtained by inspection of the
image.
\medskip

\begin{example}[Prato \cite{prato}]
\normalfont
 Let $M:={\mathbb C}^{2}\setminus
(D^{1}\times D^{1})$, where $D^{1}$ is the closed unit disc in
$\mathbb C$. If we consider $M$ as a symplectic submanifold of
$\Bbb C ^2$ with its standard symplectic structure, then we have
on $M$ we have a natural globally Hamiltonian action of $T^{2}$
given by $(e^{i{\theta}_{1}},e^{i{\theta}_{2}}) \cdot
(z_{1},z_{2}):=(e^{i{\theta}_{1}}z_{1}, e^{i{\theta}_{2}}z_{2})$.
The momentum map for this action is $(z_{1},z_{2})\mapsto
({|z_{1}|}^{2},{|z_{2}|}^{2})/2$. This momentum map has connected
fibers. Denote by $\mathbb{R}_+ := \{ x \in \mathbb{R} \mid x \geq
0 \}$. The image of the momentum map is ${\mathbb R}^{2}_{+}
\setminus \{(x,y) \mid x\leq 1/2 \,\, \text{and} \,\, y\leq 1/2
\}$ and $\mathbf{C}(J_{T^{2}}(M^{reg}))$ is the union of $\{(x,0)\mid
x>1/2\}$ and $\{(0,y)\mid y>1/2\}$ which does not disconnect any
region in the image $J _{T^{2}} (M) $. Consequently, according to
Theorem~\ref{prima}, this momentum map is open onto its image and
has a locally convex image, in agreement with Theorem~\ref{conv}.
\quad $\blacklozenge$
\end{example}

\begin{example}[Karshon and Lerman \cite{karler}]
\normalfont Let $M_{1}:=T^{2} \times U$ where $T^{2}$ is the two
dimensional torus and $U$ is the subset of ${\mathbb R}^{2}$
obtained by removing the origin and the positive $x$-axis. $M _1
$ is a symplectic manifold when viewed as an open submanifold of
the cotangent bundle of $T ^2$. The restriction to $M _1 $ of the
lifted action of $T ^2 $  on  its cotangent bundle has as momentum
map  the projection onto $U$. Hence, the image of this momentum
map is ${\mathbb R}^{2}$ minus the origin and the positive
$x$-axis.

Let $M_{2}$ be the symplectic manifold ${\mathbb C}^{2}$
minus the points whose first coordinate is nonzero. The momentum
map for the $T^{2}$ action on $M_{2}$ is given by $(z,w)\mapsto
(|z |^{2},| w |^{2})/2$ and the image is the set $\{(x,y)\in
{\mathbb R}^{2} \mid x>0, y\geq 0\}$.

Gluing these two spaces $M_1 $ and $M_2 $ along the pre-images
of the positive quadrant we obtain another globally
Hamiltonian $T^{2}$-space with a momentum map with connected
fibers whose image is
${\mathbb R}^{2}$ minus the origin. It is easy to see that
$\mathbf{C}(J_{T^{2}}(M^{reg}))$ is the open positive $x$-axis which
disconnects regions in ${\mathbb R}^{2}$.
Theorem~\ref{prima}  implies that this momentum map is not open
onto its image.
\quad $\blacklozenge$
\end{example}

To prove a converse to Proposition \ref{conectregion} in the case
when the fibers are not connected we need a few more topological
facts. A metric space is called a {\bfi generalized
continuum\/} if it is locally compact and connected. In a
topological space a {\bfi quasi-component\/} of a point is the
intersection of all closed-and-open sets that contain that point.
A topological space is called {\bfi totally disconnected\/} if the
quasi-component of any point consists of the point itself. A
continuous map $f:X \rightarrow Y$ is called {\bfi light\/} if all
fibers $f^{-1}(y)$ are totally disconnected. We say that a subset
of a topological space is {\bfi non-dense\/} if and only if it
contains no open subsets. Whyburn \cite{whyburn} (page 94) proved
the following result called ``the extension of openness''.

\begin{theorem}[Whyburn]
\label{Whyburn}
Let $X$ and $Y$ be locally connected
generalized continua and let $f:X \rightarrow Y$ be an onto light
mapping which is open on $X \backslash f^{-1}(F)$, where $F$ is a
closed non-dense set in $Y$ which separates no region in $Y$ and
is such that $f^{-1}(F)$ is non-dense. Then $f$ is open on $X$.
\end{theorem}

If the fibers of $J_T $ are not connected we still need a control
on the connected components of the fibers of $J_T $ in order to
have a similar result to that in Proposition \ref{connected
fibers} (which is a converse of Proposition \ref{conectregion}).

\begin{definition}
Let $J_{T}:M\rightarrow \frak{t}^{*}$ be the momentum map of a
torus action on a connected symplectic manifold $(M,\omega )$. We
say that $J_T $ satisfies the {\bfi connected component fiber
condition\/}
\begin{description}
\item[(CCF)]  if $J_T(x) = J_T(y) $ and $E_x \cap M^{reg} \neq
\varnothing $, then $E_y \cap M^{reg} \neq \varnothing $, where
$E_x $ and $E_y $ are the connected components of the fiber $J_T
^{-1}(J_T(x))$ that contain $x $ and $y $ respectively.
\end{description}
\end{definition}

Recall that $M_{J_T}$ is the quotient topological space whose
points are the equivalence classes given by the connected
components of the fibers of $J_T$. Denote by $\pi_{J_T}: M
\rightarrow M_{J_T}$ the canonical projection.

\begin{proposition}
\label{disconnected fibers} Suppose that $M_{J_{T}}$ is a
Hausdorff space, $J_{T}(M)$ is locally compact,
$\mathbf{C}J_{T}(M^{reg})$ does not disconnect any region in
$J_{T}(M)$, and $J_{T}$ also satisfies condition (CCF). Then
$J_{T}$ is open onto its image.
\end{proposition}

\noindent\textbf{Proof.\ \ } To prove the result we will show that
the conditions of Whyburn's Theorem are satisfied where we take
$f$ to be $\widetilde{J}_{T}:M_{J_T} \rightarrow J_T(M) \subset
\mathfrak{t}^\ast$, that is, the quotient map uniquely defined by  $J_T =
\pi_{J_T} \circ \widetilde{J}_{T}$ and $F:=\mathbf{C}J_{T}(M^{reg})$,
which is closed and non-dense in $J_T(M)$ by Lemma \ref{affine space}.

By hypothesis $M_{J_{T}}$ is a Hausdorff space. Using the fact
that $M$ is locally compact and $\pi_{J_{T}}$ is open (see Lemma
\ref{primatop} \textbf{(i)}) we obtain that $M_{J_{T}}$ is locally
compact. Since $M $ is connected, its quotient $M_{J_{T}}$ is
connected. Proposition \ref{metrica} guaranties that $M_{J_{T}}$
is a metric space since it is Hausdorff. Therefore, $M_{J_{T}}$ is
a generalized continuum. The same is true for $J_{T}(M)$. Both are
locally connected since $M$ is path connected.

Now we prove that $\widetilde{J}_{T}:M_{J_T} \rightarrow J_T(M)$
is a light map. For this, take $v \in J_T(M) $. We want to show
that $\widetilde{J}_T ^{-1}(v) $ is totally disconnected. Let $[x]
\in \widetilde{J}_T ^{-1}(v)$ be arbitrary and choose $x \in M $ a
representative of this class. Since $J_T $ has local convexity
data and is locally fiber connected (see Theorem \ref{conuri}), we
can find a small neighborhood $U_x $ of $x $ in $M$ such that
$\pi_{J_T}(U_x)$ is open in $M_{J_T}$ and is such that
$\widetilde{J}_T^{-1}(v) \cap \pi_{J_T}(U_x) = [x]$. Thus,
$\widetilde{J}_{T}:M_{J_T} \rightarrow J_T(M)$ is a light map.

Now we prove that $\widetilde{J}_{T}^{-1}(F)$ is non-dense in
$M_{J_T}$. By contradiction, suppose that this is not true. Then
there exists an open set $U \subset \widetilde{J}_{T}^{-1}(F)$.
Because $\pi_{J_T}(M^{reg})$ is dense in $M_{J_T}$ we have that $U
\cap \pi_{J_T}(M^{reg}) \neq \varnothing$. Thus there exists an
element $[x]\in \widetilde{J}_{T}^{-1}(F) \cap \pi_{J_T}(M^{reg})$
and hence there is an $x \in M^{reg} $ such that $J_T(x) =
\widetilde{J}_T([x]) \in F $. This contradicts the definition of
$F$.

The last thing to be verified is that $\widetilde{J}_{T}$
restricted to  $M_{J_{T}} \backslash \widetilde{J}_{T}^{-1}(F)$ is
an open map. To see this note that the inclusion
\begin{equation*}
\widetilde{J}_{T}^{-1}(F) =
\widetilde{J}_{T}^{-1}(\mathbf{C}J_{T}(M^{reg}))
=\mathbf{C}\widetilde{J}_{T}^{-1}
(J_{T}(M^{reg}))\subset\mathbf{C}\pi_{J_{T}}(M^{reg})
\end{equation*}
shows that $M_{J_{T}} \backslash \widetilde{J}_{T}^{-1}(F)\supset
\pi_{J_{T}}(M^{reg})$. Now we shall prove the reverse inclusion.
Let $[x] \in M_{J_{T}} \backslash \widetilde{J}_{T}^{-1}(F)$. If
the connected component $E_x $ of the fiber $J_T ^{-1}(J_T(x))$
intersects $M^{reg}$, then $[x] \in \pi_{J_T}(M^{reg})$. If not,
then we have $J_T(x) = \widetilde{J}_T([x]) \notin F $, that is,
$J_T(x) \in \mathbf{C}F = J_T(M^{reg})$. Therefore, there is some
$y \in M^{reg} $ such that $J_T(x) = J_T(y)$. By condition (CCF),
$E_x \cap M^{reg} \neq \varnothing $ and hence $[x] \in
\pi_{J_T}(M^{reg})$ which proves the equality $M_{J_{T}}
\backslash \widetilde{J}_{T}^{-1}(F) = \pi_{J_{T}}(M^{reg})$.
Thus, since $\widetilde{J}_{T}$ is open in $\pi_{J_{T}}(M^{reg})$
the last requirement of Whyburn's Theorem is verified. \quad
$\blacksquare$

\medskip

To summarize, in the situation for non-connected fibers we obtain
the following result.

\begin{theorem}
\label{adoua} Let $J_{T}:M\rightarrow \frak{t}^{*}$ be the
momentum map of a torus action on a connected symplectic manifold
$(M,\omega )$. Suppose that $M_{J_{T}}$ is a Hausdorff space. 

Then
$J_{T}$ is open onto its image if and only if $J_{T}(M)$ is
locally compact,  $\mathbf{C}J_{T}(M^{reg})$ does not disconnect
any region in $J_{T}(M)$, and $J_T $ satisfies condition (CCF).
Moreover, under these hypotheses, the image of the momentum map is
locally convex and locally polyhedral.
\end{theorem}

\noindent\textbf{Proof.\ \ } The only thing that remains to be
shown is that if $J_T: M \rightarrow J_T(M) $ is open onto its
image then condition (CCF) holds. Suppose that the condition (CCF)
does not hold. So there exists a fiber with at least two connected
components $E_x$ and $E_y$ such that $E_x \cap M^{reg} \neq
\varnothing $ and $E_y \cap M^{reg}= \varnothing $. 

Consequently,
we can suppose that $x\in M^{reg}$ and $y$ is contained in a lower
stratum of the $T$ action. Then we have the strict inclusion
$C_{J_{T}(y)}=v+\mathfrak{t}_{0}^{\perp}+
\operatorname{cone}(\mathcal{P}_{V})\subset C_{J_{T}(x)}=v+
\mathfrak{t}_{1}^{\ast}$, where $v=J_{T}(x)=J_{T}(y)$ (see Theorem
\ref{conuri} and the proof of Lemma \ref{affine space}). By condition (VN), there exist open 
neighborhoods $U _x $ and $U _y $
of $x $ and $y$, respectively, such that $J _T (U _x )$  is an
open ball in $\mathfrak{t} _1 ^\ast  $ centered at  $v $  and  $J
_T (U _y )$ is the intersection of an open ball in $\mathfrak{t}
_1 ^\ast  $ centered at  $v $ with a closed proper cone in
$\mathfrak{t} _1 ^\ast$ with vertex $v$. This contradicts the
openness onto its image of $J_{T}$. \quad $\blacksquare$

\medskip

In the case of a compact, connected, and  non-Abelian group, the momentum map
$J_{G}: M\rightarrow \frak{g} ^{\ast}$ is, in general, not open
onto its image even if it is a proper map. Nevertheless, we will
show that similar results to those obtained in the Abelian case
hold for the quotient map $j_{G}:M\rightarrow \frak{g}^{\ast}/G
\simeq {\mathfrak{t}^{\ast}_+}$, where $j_{G}=\pi_{G}\circ J_{G}$
and $\pi_{G}:\frak{g}^{\ast}\rightarrow\mathfrak{t}^{\ast}_+$ is
the projection map which is always proper if $G$ is compact.

The quotient map $j_G $ has local convexity data due to the
following result of Sjamaar \cite[Theorem 6.5]{sjamaar}.

\begin{theorem}[Sjamaar]
Let $M$ be a connected Hamiltonian $G$-manifold. Then for every
$x\in M$ there exist a unique, closed, polyhedral convex cone
$C_{x}$ in $\mathfrak{t}^\ast _+$ with vertex at $j_{G}(x)$ such
that for every sufficiently small $G$-invariant neighborhood $U$
of $x$ the set $j_{G}(U)$ is an open neighborhood of $j_{G}(x)$ in
$C_{x}$.
\end{theorem}

Using Lerman's \cite{lerman} symplectic cut technique, Knop
\cite[Theorem 5.1]{knop} proved that $j_G $ is locally fiber
connected. As a consequence of these results we have the following
theorem.

\begin{theorem}[Knop-Sjamaar] The map $j_{G}$ is locally fiber
connected and has local convexity data.
\end{theorem}

Using the above theorem and Theorem~\ref{localglobal} we obtain
the following generalization of Kirwan's convexity result
\cite{kirwan convexity}. In the next statement we will use the map
$ \widetilde{j}_G:M/G \rightarrow {\mathfrak{t}^{\ast}_+}$
defined by the identity $j_G=  \widetilde{j}_G \circ \pi$, where
$\pi:M \rightarrow  M/G $ is the projection. We will say that the
$G$-equivariant  momentum map $J_G:M \rightarrow \mathfrak{g}^\ast
$ is {\bfi  $G$-open } onto its image whenever $\widetilde{j}_G $
is open onto its image.

\begin{theorem}
Let $M$ be a connected Hamiltonian $G$-manifold with
$G$ a compact connected Lie group. If the momentum map $J_{G}$ is closed
then $J_{G}(M)\cap\mathfrak{t}^{\ast}_+$ is a closed convex
locally polyhedral set. Moreover, $J _G $ is $G$-open onto its
image and all its fibers are connected.
\end{theorem}

\noindent\textbf{Proof.\ \ } As a direct application of Theorem
~\ref{localglobal} we obtain that $j_{G}$ is open onto its image
and consequently that $J_{G}$ is $G$-open  onto its image. Additionally, the set
$j_{G}(M)=J_{G}(M)\cap\mathfrak{t}^{\ast}_+$ is a closed convex
locally polyhedral set and $j_{G}$ has connected fibers. 

It
remains to be proved that $J_{G}$ has connected fibers. To see this,  note that since
$j_{G}$ has connected fibers, the pre-images
$J_{G}^{-1}(\mathcal {O}_{\mu})$ are connected as topological subspaces of $M$ for
every coadjoint orbit
$\mathcal{O}_{\mu}\subset J_{G}(M)$.
 Note now that $J _G^{-1}(\mathcal{O}_{\mu})$
can also be endowed with the initial topology
induced by the map 
\[
\begin{array}{cccc}
J _G^{\mathcal{O}_{\mu}} :& J _G^{-1}(\mathcal{O}_{\mu}) & \longrightarrow &
\mathcal{O}_{\mu}\\
	&z &\longmapsto & J _G(z),
\end{array}
\]
where the orbit $\mathcal{O}_{\mu} $ comes with its orbit smooth structure
induced by the homogeneous  manifold $G/G _{\mu} $. Since $G$ is compact, the
orbit 
$\mathcal{O}_{\mu} $ is an embedded submanifold of $\mathfrak{g}^\ast$ and hence 
the initial topology for $J _G (\mathcal{O}_{\mu}) $ is weaker than the subspace
topology. Indeed, the sets of the form $(J _G^{\mathcal{O}_{\mu}}) ^{-1} (U\cap 
\mathcal{O}_{\mu})=J _G^{-1}(U)\cap J_G ^{-1}(\mathcal{O}_{\mu})$, with $U$ open
in $\mathfrak{g}^\ast$, form a subbasis of the initial topology of  $J
_G^{\mathcal{O}_{\mu}} $ and since they are open in $M$ by the continuity of $J _G
$, the claim follows. Therefore, $J _G^{-1}(\mathcal{O}_{\mu}) $ is also connected
for the initial topology. Proposition 8.4.1 in~\cite{OrRa2004} states that 
if  $ J _G^{-1}(\mathcal{O}_{\mu})  $ is endowed with its initial topology then 
the map 
\[
\begin{array}{cccc}
f : & G \times  _{G_{\mu} } J _G^{-1}(\mu) & \longrightarrow &
J _G^{-1}(\mathcal{O}_{\mu})\\
  &[g,z] &\longmapsto & g \cdot  z, 
\end{array}
\] 
is a homeomorphism, where $G \times  _{G_{\mu} } J _G^{-1}(\mu)$ denotes the 
orbit space of the free and continuous action $h \cdot (g, z):= (gh, h
^{-1}\cdot  z )$, $h \in  G _\mu $, $g \in  G $,  $z \in   J _G ^{-1}(\mu) $, of
the compact connected group (see Theorem 3.3.1 in~\cite{dk})
$G _\mu $  on the product $G \times J _G ^{-1}(\mu) $.
The set $J _G ^{-1}(\mu) $ is considered with its subspace topology. Let $\pi
_\mu: G \times  J _G ^{-1}(\mu) \rightarrow G \times  _{G _\mu} J _G ^{-1}(\mu) $ 
be the continuous and open projection. Since the fibers of $\pi_\mu $ are
connected (they are homeomorphic to $G _\mu$) and $\pi _\mu  $  is open it follows
that the pre-image of any connected set is connected by Lemma~\ref{conex}.
Therefore $G \times  J _G ^{-1}(\mu) $ is connected and hence so is $J _G
^{-1}(\mu) $ since $G$ is connected.\quad $\blacksquare$

\begin{remark}
\normalfont
We emphasize that $J _G  $ is $G$-open but not open in general.
See~\cite{montaldi}  for a counterexample.
\end{remark}

Analyzing the proofs of the results leading to Theorem \ref{prima}
and using the natural definition of the (CCF) condition for the
map $j_{G}$ we have the following two results (see Remark
\ref{topological convexity remark}) which are the non-Abelian
analogs of Theorem~\ref{prima} and \ref{adoua}.

\begin{theorem}
\label{atreia}
Let $G$ be a compact connected Lie group and $M$ be a connected
Hamiltonian $G$-manifold with equivariant momentum map
$J_{G}:M\rightarrow \frak{g} ^{\ast}$. Suppose that $J_{G}$
has connected fibers. Then $J_{G}$ is $G$-open onto its image if
and only if $\mathbf{C}((\pi_{G}\circ J_{G})(M^{reg}))$ does not
disconnect any region in $J_{G}(M)\cap \mathfrak{t}_{+}^{\ast }$.
Moreover, in this context, the image $J_{G}(M)\cap
\mathfrak{t}_{+}^{\ast }$ is a locally convex and locally
polyhedral set.
\end{theorem}

\begin{theorem}
\label{apatra} Let $G$ be a compact connected group and $M$ be a
connected Hamiltonian $G$-manifold with the momentum map
$J_{G}:M\rightarrow \frak{g} ^{\ast}$. Suppose that
$(M/G)_{j_{G}}$ is a Hausdorff space. Then $J_{G}$ is $G$-open
onto its image if and only if $J_{G}(M)$ is locally compact,
$\mathbf{C}(\pi_{G}\circ J_{G}(M^{reg}))$ does not disconnect any region in
$J_{G}(M)\cap \mathfrak{t}_{+}^{\ast }$, and $j_{G}$ satisfies
the (CCF) condition. Moreover, in this context, the image $J_{G}(M)\cap
\mathfrak{t}_{+}^{\ast }$ is a locally convex and locally
polyhedral set.
\end{theorem}

\section{Convexity for Poisson actions of compact Lie groups}
\label{Convexity for Poisson actions of compact Lie groups}

In this section we shall prove a generalization of the Flaschka-Ratiu
convexity theorem (Theorem 4.39 in~\cite{flaschkaratiu}) which has as
one of its main consequences the convexity theorem for Poisson actions
of compact connected Poisson Lie groups on compact connected symplectic
manifolds. We hasten to add that the hypotheses of this theorem do
\textit{not\/} imply that the compact Lie group action is necessarily a
Poisson action of a Poisson Lie group and that this theorem can be applied
in other  situations. The generalization below will require only that
a certain map is closed (not even properness is assumed) instead of assuming
compactness of the symplectic manifold.

Let $K$ be a compact connected semisimple Lie group. Since any
compact Lie group $K$ is the commuting product $(Z_K)_0 K_{ss}$ of the connected
component of the identity $(Z_K)_0 $ of the center $Z_K$ and of a closed semisimple
subgroup $K_{ss}$ (see Theorem 4.29 in~\cite{knapp}) we will provide in what
follows, without any loss of generality, the proofs for the compact semisimple case,
even though the results will be stated for general compact Lie groups.
Since $K$ admits
a complexification, we can think of it as the compact real form of a
connected complex semisimple Lie group $G$. Denote by
$G^\mathbb{R}$ the real Lie group underlying $G$ and let
$G^\mathbb{R}=KAN$ be its Iwasawa decomposition. Denote by
$\mathfrak{g}^ \mathbb{R}, \mathfrak{k}, \mathfrak{a}$, and
$\mathfrak{n}$ the real Lie algebras of $G^\mathbb{R}, K, A $, and
$N $, respectively. Then $\mathfrak{g}^\mathbb{R} =
\mathfrak{k}\oplus \mathfrak{a} \oplus \mathfrak{n}$ is the
Iwasawa decomposition of $\mathfrak{g}^ \mathbb{R}$.  If
$\mathfrak{t} = i\mathfrak{a}$ then $T =\exp\mathfrak{t}$ is a
maximal torus of $K$. Define $B : = AN$ whose Lie algebra is
$\mathfrak{b}:= \mathfrak{a} \oplus \mathfrak{n}$.

Let $\kappa$ be the Killing form of $\mathfrak{g}$. Its imaginary
part $\operatorname{Im} \kappa$ is a non-degenerate
invariant symmetric bilinear form on $\mathfrak{g}^ \mathbb{R}$.
Since $\operatorname{Im} \kappa( \mathfrak{k}, \mathfrak{k}) =
\operatorname{Im}(\mathfrak{b}, \mathfrak{b}) = 0$, the vector
spaces $\mathfrak{k}$ and $\mathfrak{b}$ are dual to each other
relative to $\langle\,, \rangle : =\operatorname{Im} \kappa$. The
Cartan decomposition $G^\mathbb{R} = PK$ defines the Cartan
involution $\tau: G^\mathbb{R} \rightarrow G^\mathbb{R}$. Define
$g ^\ast: = \tau( g ^{-1}) $ for any $g \in G $. The derivative of
these maps at the identity will be denoted by the same symbols. The
map $\tau: \mathfrak{g}^ \mathbb{R} \rightarrow \mathfrak{g}^
\mathbb{R}$ has eigenvalues $\pm 1$. The $+1 $ eigenspace is
$\mathfrak{k}$ and the $-1 $ is denoted by $\mathfrak{p}$. Note
that $\left\langle \,, \right\rangle$ identifies
$\mathfrak{k}^\ast$ with $\mathfrak{p}$. The exponential map is a
diffeomorphism from $\mathfrak{p}$ to $P $.

Let $\mathfrak{a}_+$ be the positive Weyl chamber in
$\mathfrak{a}\cong\mathfrak{t}^\ast$ corresponding to the subgroup
$B$, let $\mathfrak{a}_0$ be the interior of $\mathfrak{a}_+$, and
set $A_+=\exp\mathfrak{a}_+$ and
$A_0=\exp\mathfrak{a}_0$. Let $W^1, W^2, ...$ be all the closed
walls of varying dimensions of the positive Weyl chamber
$\mathfrak{a}_+$; there are only finitely many such walls. If
$P^i$ is the subspace of $\mathfrak{a}$ spanned by $W^i$, then
$W^i$ is closed in $P^i$ and $W^i_0$ denotes the interior of
$W^i$ in $P^i$. We will also use the notation $W^0 =
\mathfrak{a}_+$ and $W^0_0 = \mathfrak{a}_0$; note that $P^0 =
\mathfrak{a}$.

\begin{theorem} Suppose that a compact connected Lie
group $K$ acts on a paracompact connected symplectic manifold
$(M, \omega)$ and that a maximal torus $T$ of $K$ acts in a
Hamiltonian fashion with equivariant momentum map $J_T: M
\rightarrow \mathfrak{t}^\ast$.  Suppose there exists a closed map
$\mathcal{P}:M\rightarrow\mathfrak{p}$ with the following
properties:
\begin{description}
\item[(i)] $\mathcal{P}$ is equivariant {\rm (}with respect to the
adjoint action of $K$ on $\mathfrak{p}${\rm )};
\smallskip
\item[(ii)] for every $x\in M$, $T_x\mathcal{P}(T_xM)=
\mathfrak{k}_x^{\operatorname{ann}} : = \{ \mu \in \mathfrak{p}
\mid \left\langle \mu, \xi \right\rangle = 0 \text{~for all~} \xi
\in \mathfrak{k}\}$;
\smallskip
\item[(iii)] for every $x\in M$, the kernel of $T_x\mathcal{P}$ is
\[
(\mathfrak{k} \cdot x)^\omega : = \{v \in T_x M \mid \omega(x)(v,
\xi_M(x))= 0 \text{ for all } \xi \in \mathfrak{k}\};
\]
\item[(iv)] the restriction of $\mathcal{P}$ to
$\mathcal{P}^{-1}(\mathfrak{a}_+)$ is proportional to $J_T$.
\end{description}
Then $\mathcal{P}(M)\cap\mathfrak{a}_+$ is a closed convex set. If
$M$ is compact, the set $\mathcal{P}(M)\cap\mathfrak{a}_+$ is a
compact convex polytope.
\end{theorem}

\noindent\textbf{Proof.\ \ } To prove the statement we shall use
the technique of symplectic cross sections for the map $\mathcal{P}$
(see \cite{convexity, dazord, hilnebplank}).  Two situations may arise: if the
intersection $\mathcal{P}(M)\cap \mathfrak{a} _0 $ is nonempty then the
symplectic cross section can be taken to be $Y:= \mathcal{P} ^{-1}
(\mathfrak{a}_0) $. Otherwise, $\mathcal{P}(M)\cap \mathfrak{a} _0 =\varnothing$ and
$Y:= \mathcal{P} ^{-1} (W _0 ^i)$ is a symplectic cross section, where $W ^i $ is
the unique closed wall such that  $\mathcal{P} (M) \cap  \mathfrak{a}_+ \subset W ^i
$ and  $\mathcal{P} (M)\cap W _0^i\neq \varnothing $.

The hypotheses in the theorem ensure via Proposition 4.27 in~\cite{flaschkaratiu}
that the set $Y$ defined above is a symplectic cross section and that, additionally,
the closure of $\mathcal{P} (Y) $ in $\mathfrak{a}_+ $ is $ \mathcal{P}(M)\cap
\mathfrak{a}_+ $. Consequently, it suffices to show that
$\mathcal{P} (Y) $ is convex in order to prove the convexity of  $\mathcal{P}(M)
\cap  \mathfrak{a} _+ $.

The  proportionality of $\mathcal{P} |_Y $ to  $ J_T $ guarantees
that $\mathcal{P}|_Y $ is locally fiber connected and has local
convexity data. Notice that $\mathcal{P}|_Y:Y \rightarrow
W ^i_0  $ is a closed map. Then, using Remark~\ref{restriction} it
follows that $\mathcal{P}(Y) $ is convex, closed, and locally polyhedral in $W _0 ^i
$. 

If $M $ is compact then Proposition 4.38 in~\cite{flaschkaratiu} states that
$\overline{\mathcal{P}(M)}$ is polyhedral and that  it has only a finite number of
faces. Thus $\mathcal{P}(M)\cap\mathfrak{a}_+$ is a compact convex polytope.
\quad $\blacksquare$

\medskip

\noindent\textbf{Acknowledgments}   
This research was partially supported by the
European Commission through funding for the Research Training Network
\emph{Mechanics and Symmetry in Europe} (MASIE) as well as the Swiss National
Science Foundation. The  authors are thankful for the excellent
working conditions provided by the Bernoulli Center where  part of this work was
done during the program ``Geometric Mechanics and Its Applications". Petre Birtea
has been supported by a grant of the R\'egion de Franche-Comt\'e (Convention
051004-02) during his stay at the Universit\'e de Franche-Comt\'e, Besan\c{c}on,
which made possible part of this collaboration.

\noindent {\sc P. Birtea} \\
Departamentul de Matematic\u a, Universitatea de Vest,
RO--1900 Timi\c soara, Romania.\\
Email: {\sf birtea@math.uvt.ro}
\medskip

\noindent {\sc J.-P. Ortega} \\
D\'epartement de Math\'ematiques de Besan\c{c}on,
Universit\'e de Franche-Comt\'e, UFR des Sciences et
Techniques, 16 route de Gray, F--25030 Besan\c{c}on c\'edex,
France.\\
Email: \texttt{Juan-Pablo.Ortega@math.univ-fcomte.fr}
\medskip

\noindent {\sc T.S. Ratiu}\\
Centre Bernoulli,
{\'E}cole Polytechnique F{\'e}d{\'e}rale de Lausanne,
CH--1015 Lausanne,
Switzerland.\\
Email: {\sf tudor.ratiu@epfl.ch}
\end{document}